\documentclass[11pt]{article}
\usepackage{amssymb}
\usepackage{amsmath}
\usepackage[all]{xy}
\usepackage{times}
\usepackage{graphicx}

\title{Descriptive complexity of countable unions of Borel rectangles\indent}
\author{Dominique LECOMTE and Miroslav ZELENY$^1$}
\date{\today}

\def\ufootnote#1{\let\savedthfn\thefootnote\let\thefootnote\relax
\footnote{#1}\let\thefootnote\savedthfn\addtocounter{footnote}{-1}}

\newcommand{\Ana}{{\it\Sigma}^{1}_{1}}

\newcommand{\borel}{{\bf\Delta}^{1}_{1}}

\newcommand{\boraone}{{\bf\Sigma}^{0}_{1}}
\newcommand{\boratwo}{{\bf\Sigma}^{0}_{2}}
\newcommand{\borathree}{{\bf\Sigma}^{0}_{3}}

\newcommand{\boraxi}{{\bf\Sigma}^{0}_{\xi}}

\newcommand{\borone}{{\bf\Delta}^{0}_{1}}
\newcommand{\bortwo}{{\bf\Delta}^{0}_{2}}
\newcommand{\borthree}{{\bf\Delta}^{0}_{3}}

\newcommand{\bormone}{{\bf\Pi}^{0}_{1}}

\newcommand{\bormtwo}{{\bf\Pi}^{0}_{2}}

\newcommand{\boraxp}{{\bf\Sigma}^{0}_{\xi +1}}

\newcommand{\bormlxi}{{\bf\Pi}^{0}_{<\xi}}

\newcommand{\borath}{{\bf\Sigma}^{0}_{\theta}}

\newcommand{\bormxi}{{\bf\Pi}^{0}_{\xi}}

\newcommand{\bormen}{{\bf\Pi}^{0}_{\eta_n}}
\newcommand{\bormem}{{\bf\Pi}^{0}_{\eta_m}}

\newcommand{\borxi}{{\bf\Delta}^{0}_{\xi}}

\newcommand{\bormth}{{\bf\Pi}^{0}_{\theta}}

\newtheorem{thm} {Theorem} [section]
\newtheorem{cor} [thm] {Corollary}
\newtheorem{lem} [thm] {Lemma}
\newtheorem{prop} [thm] {Proposition}
\newtheorem{defi} [thm] {Definition}

\def\Baire{\omega^\omega}

\begin{document}

\maketitle

\centerline{$\bullet$ Universit\' e Paris 6, Institut de Math\'ematiques de Jussieu, Projet Analyse Fonctionnelle}

\centerline{Couloir 16-26, 4\`eme \'etage, Case 247, 4, place Jussieu, 75 252 Paris Cedex 05, France}

\centerline{dominique.lecomte@upmc.fr}\bigskip

\centerline{$\bullet$ Universit\'e de Picardie, I.U.T. de l'Oise, site de Creil,}

\centerline{13, all\'ee de la fa\"\i encerie, 60 107 Creil, France}\bigskip

\centerline{$\bullet^1$ Charles University, Faculty of Mathematics and Physics, Department of Mathematical Analysis}

\centerline{Sokolovsk\'a 83, 186 75 Prague, Czech Republic}

\centerline{zeleny@karlin.mff.cuni.cz}\bigskip\bigskip\bigskip\bigskip\bigskip\bigskip

\ufootnote{{\it 2010 Mathematics Subject Classification.}~Primary: 03E15, Secondary: 54H05}

\ufootnote{{\it Keywords and phrases.}~Borel chromatic number, Borel class, coloring,  product, rectangle}

\ufootnote{{\it Acknowledgements.}~The results were partly obtained during the second author's 
stay at the Universit\'e Paris 6 in May 2013. The second author thanks the Universit\'e Paris 6 for the hospitality.}

\ufootnote{The research was supported by the grant GA{\v C}R P201/12/0436 for the second author.}

\noindent {\bf Abstract.} We give, for each countable ordinal $\xi\!\geq\! 1$, an example of a 
$\bortwo$ countable union of Borel rectangles that cannot be decomposed into countably many 
$\bormxi$ rectangles. In fact, we provide a graph of a partial injection with disjoint domain and range, which is a difference of two closed sets, and which has no $\borxi$-measurable countable coloring. 

\vfill\eject

\section{$\!\!\!\!\!\!$ Introduction}\indent

 In this paper, we work in products of two Polish spaces. One of our goals is to give an answer to the following simple question. Assume that a countable union of Borel rectangles has low Borel rank. Is there a decomposition of this union into countably many rectangles of low Borel rank? In other words, is there a map 
$r\! :\!\omega_1\!\setminus\!\{ 0\}\!\rightarrow\!\omega_1\!\setminus\!\{ 0\}$ such that 
$\bormxi\cap (\borel\!\times\!\borel )_\sigma\!\subseteq\! 
({\bf\Pi}^0_{r(\xi )}\!\times\! {\bf\Pi}^0_{r(\xi )})_\sigma$ for each $\xi\!\in\!\omega_1\!\setminus\!\{ 0\}$?\bigskip

 By Theorem 3.6 in [Lo], a Borel set with open vertical sections is of the form 
$(\borel\!\times\!\boraone )_\sigma$. This leads to a similar problem: is there a map 
$s\! :\!\omega_1\!\setminus\!\{ 0\}\!\rightarrow\!\omega_1\!\setminus\!\{ 0\}$ such that, for each 
$\xi\!\in\!\omega_1\!\setminus\!\{ 0\}$, $\bormxi\cap (\borel\!\times\!\boraone )_\sigma\!\subseteq\! 
({\bf\Pi}^0_{s(\xi )}\!\times\! {\bf\Sigma}^0_1)_\sigma$?\bigskip

 The answer to these questions is negative: 
 
\begin{thm} Let $1 \leq \xi < \omega_1$. Then there exists a partial map 
$f\! :\!\omega^\omega\!\rightarrow\!\omega^\omega$ such that the complement 
$\neg\mbox{Gr}(f)$ of the graph of $f$ is $\bormtwo$ but not $(\boraxi\!\times\! \borel )_{\sigma}$.\end{thm}
 
 In fact, we prove a result related to $\borxi$-measurable countable colorings. A study of such colorings is made in [L-Z]. It was motivated by the $\mathbb{G}_0$-dichotomy (see Theorem 6.3 in [K-S-T]). More precisely, let $B$ be a Borel binary relation having a Borel countable coloring (i.e., a Borel map $c\! :\! X\!\rightarrow\!\omega$ such that $c(x)\!\not=\! c(y)$ if $(x,y)\!\in\! B$). Is there a relation between the Borel class of $B$ and that of the coloring?  In other words, is there a map 
$k\! :\!\omega_1\!\setminus\!\{ 0\}\!\rightarrow\!\omega_1\!\setminus\!\{ 0\}$ such that any $\bormxi$ binary relation having a Borel countable coloring has in fact a ${\bf\Delta}^0_{k(\xi )}$-measurable countable coloring, for each $\xi\!\in\!\omega_1\!\setminus\!\{ 0\}$? Here again, the answer is negative:

\begin{thm} Let $1 \leq \xi < \omega_1$. Then there exists a partial injection with disjoint domain and range $i\! :\!\omega^\omega\!\rightarrow\!\omega^\omega$ whose graph is the difference of two closed sets, and has no $\borxi$-measurable countable coloring.\end{thm}

 These two results are consequences of Theorem 4 in [M\'a] and its proof. This latter can also be used positively, to produce examples of graphs of fixed point free partial injections having reasonable chances to characterize the analytic binary relations without $\borxi$-measurable countable coloring. We will see in Section 4 that such a characterization indeed holds when 
$\xi\! =\! 3$, and give an example much simpler than the one in [L-Z]. In Section 2, we give a proof of Theorem 4 in [M\'a], in $\omega^\omega$ instead of $2^\omega$, and also prove some additional properties needed for the construction of our partial maps. In Section 3, we prove Theorems 1.1 and 1.2. At the end of Section 4, we show that Theorem 1.2 is optimal in terms of descriptive complexity of the graph, and also give a positive result concerning the first two problems in the case of finite unions of rectangles.

\section{$\!\!\!\!\!\!$ M\'atrai sets}\indent

 Before proving our version of Theorem 4 in [M\'a], we need some notation, definition, and a few basic facts. The maps with closed graph will be of particular interest for us.
 
\begin{lem} Let $(X_i)_{i\in\omega}$, $(Y_i)_{i\in\omega}$ be sequences of metrizable spaces, and, for each $i\!\in\!\omega$, $f_i\! :\! X_i\!\rightarrow\! Y_i$ be a partial map whose graph is a closed subset of $X_i\!\times\! Y_i$. Then the graph of the partial map 
$f\! :=\!\Pi_{i\in\omega}~f_i\! :\!\Pi_{i\in\omega}~X_i\!\rightarrow\!\Pi_{i\in\omega}~Y_i$ is closed.
\end{lem}

\noindent\bf Proof.\rm\ Let $(x^j)_{j\in\omega}$ be a sequence of elements of 
$\Pi_{i\in\omega}~X_i$ converging to $x\! :=\! (x_i)_{i\in\omega}$ such that  
$\big( f(x^j)\big)_{j\in\omega}$ converges to $y\! :=\! (y_i)_{i\in\omega}\!\in\!\Pi_{i\in\omega}~Y_i$. Then $y_i\! =\! f_i(x_i)$, since $\mbox{Gr}(f_i)$ is closed, for each $i\!\in\!\omega$. This implies that 
$y\! =\! f(x)$ and the proof is finished.\hfill{$\square$}\bigskip
 
\noindent\bf Notation.\rm\ Let $X$ be a set and $\cal F$ be a family of subsets of $X$. Then the symbol $\langle\cal F\rangle$ denotes the smallest topology on $X$ containing $\cal F$.\bigskip

 The next two lemmas can be found in [K] (see Lemmas 13.2 and 13.3).

\begin{lem} Let $(X,\sigma)$ be a Polish space and $F$ be a $\sigma$-closed subset of $X$. Then the topology $\sigma_F\! :=\!\langle \sigma \cup \{ F\}\rangle$ is Polish and 
$F$ is $\sigma_F$-clopen.\end{lem}

\begin{lem} Let $(\sigma_n)_{n\in\omega}$ be a sequence of Polish topologies on $X$. Then the topology $\langle \bigcup_{n \in \omega}~\sigma_n\rangle$ is Polish.\end{lem}

\begin{lem} Let $(H_n)_{n \in \omega}$ be a disjoint family of sets in a zero-dimensional Polish space $(X,\sigma)$ and $(\sigma_n)_{n\in\omega}$ be a sequence of topologies on $X$ such that\smallskip

$\sigma_0 = \sigma$, $H_0$ is $\sigma_0$-closed,\smallskip

$\sigma_{n+1}\! =\!\langle\sigma_n\cup\{ H_n\}\rangle$, $H_{n+1}$ is $\sigma_{n+1}$-closed for every $n\!\in\!\omega$.\smallskip

\noindent Then the topology $\sigma_{\infty}\! =\!\langle\bigcup_{n\in\omega}~\sigma_n\rangle$ satisfies the following properties:\smallskip

(a) $\sigma_\infty$ is zero-dimensional Polish,\smallskip

(b) ${\sigma_\infty}_{\vert X\setminus\bigcup_{n\in\omega}H_n}\! =\!
\sigma_{\vert X\setminus\bigcup_{n\in\omega}H_n}$,\smallskip

\noindent and, for every $n\!\in\!\omega$,\smallskip

(c) ${\sigma_\infty}_{\vert H_n}\! =\!\sigma_{\vert H_n}$,\smallskip

(d) $H_n$ is $\sigma_\infty$-clopen.\end{lem}

\noindent\bf Proof.\rm\ Using Lemma 2.2 we see that each topology $\sigma_n$ is Polish. Then the topology $\sigma_{\infty}$ is Polish by Lemma 2.3. Now observe that the following claim holds.\bigskip

\noindent\bf Claim.\it\ A set $G\!\subseteq\! X$ is $\sigma_\infty$-open if and only if $G$ can be written as $G = G' \cup (\bigcup_{n \in \omega} G_n \cap H_n)$, where $G', G_n$ are $\sigma$-open.\rm\bigskip

 Note that $H_n\!\in\!\boraone (\sigma_{n+1})\!\subseteq\!\boraone (\sigma_\infty )$ and 
$H_n\!\in\!\bormone (\sigma_n)\!\subseteq\!\bormone (\sigma_\infty )$, thus 
$H_n$ is $\sigma_\infty$-clopen. Thus (d) is satisfied. Let $\cal B$ be a basis for $\sigma$ made  of $\sigma$-clopen sets. Then the family 
$${\cal B}\cup\{ G \cap H_n\mid G\!\in\! {\cal B}\ \wedge\ n\!\in\!\omega\}$$ 
is made of $\sigma_\infty$-clopen sets and form a basis for $\sigma_\infty$ by the claim. This gives (a).\bigskip

 Let $G\!\in\!\boraone (\sigma_\infty )$. By the claim, we find $\sigma$-open sets $G', G_n$ such that $G\! =\! G'\cup (\bigcup_{n \in \omega}~G_n\cap H_n)$. Then 
$G\cap (X\!\setminus\!\bigcup_{n\in\omega}~H_n)\! =\! 
G'\cap (X\!\setminus\!\bigcup_{n\in\omega}~H_n)$. This implies (b). Moreover, 
$G\cap H_n\! =\! G_n\cap H_n$, and (c) holds.\hfill{$\square$}

\vfill\eject

\noindent\bf Notation.\rm\ The symbol $\tau$ denotes the product topology on $\omega^\omega$. 

\begin{defi} We say that a partial map $f\! :\!\omega^\omega\!\rightarrow\!\omega^\omega$ is 
\bf nice\rm\ if $\mbox{Gr}(f)$ is a $(\tau\!\times\!\tau )$-closed subset of 
$\omega^\omega\!\times\!\omega^\omega$.\end{defi}

 The construction of $P_\xi$ and $\tau_\xi$, and the verification of the properties 
(1)$_\xi$-(3)$_\xi$ from the next lemma, can be found in [M\'a], up to minor modifications.
 
\begin{lem} Let $1 \leq \xi < \omega_1$. Then there are $P_\xi\!\subseteq\!\omega^\omega$, and a topology $\tau_\xi$ on $\omega^\omega$ such that\smallskip

(1)$_\xi$ $\tau_\xi$ is zero-dimensional perfect Polish and 
$\tau\!\subseteq\!\tau_\xi\!\subseteq\!\boraxi (\tau)$,\smallskip

(2)$_\xi$ $P_\xi$ is a nonempty $\tau_\xi$-closed nowhere dense set,\smallskip

(3)$_\xi$ if $S\!\in\!\boraxi (\Baire ,\tau)$ is $\tau_\xi$-nonmeager in $P_\xi$, then $S$ is 
$\tau_\xi$-nonmeager in $\omega^\omega$,\smallskip

(4)$_{\xi}$ if $U$ is a nonempty ${\tau_\xi}_{\vert P_\xi}$-open subset of $P_\xi$, then we can find a $\tau_\xi$-dense $G_\delta$ subset $G$ of $U$, and a nice $(\tau_\xi ,\tau )$-homeomorphism 
$\varphi_{\xi ,G}$ from $G$ onto $\omega^\omega$,\smallskip

(5)$_{\xi}$ if $V$ is a nonempty $\tau_\xi$-open subset of $\omega^\omega$, then we can find a 
$\tau_\xi$-dense $G_\delta$ subset $H$ of $V$, and a nice $(\tau_\xi ,\tau )$-homeomorphism 
$\psi_{\xi ,H}$ from  $H$ onto $\omega^\omega$,\smallskip

(6)$_{\xi}$ if $U$ is a nonempty ${\tau_\xi}_{\vert P_\xi}$-open subset of $P_\xi$ and $W$ is a nonempty open subset of $\Baire$, then we can find a $\tau_\xi$-dense $G_\delta$ subset $G$ of $U$, a $\tau_\xi$-dense $G_\delta$ subset $K$ of $W\!\setminus\! P_\xi$, and a nice 
$(\tau_\xi ,\tau_\xi )$-homeomorphism $\varphi_{\xi ,G,K}$ from $G$ onto $K$,\smallskip

(7)$_{\xi}$ if $V,W$ are nonempty $\tau_\xi$-open subsets of $\omega^\omega$, then we can find a $\tau_\xi$-dense $G_\delta$ subset $H$ of $V\!\setminus\! P_\xi$, a $\tau_\xi$-dense $G_\delta$ subset $L$ of $W\!\setminus\! P_\xi$, and a nice $(\tau_\xi ,\tau_\xi )$-homeomorphism $\psi_{\xi ,H,L}$ from $H$ onto $L$.\end{lem}

\noindent\bf Proof.\rm\ We proceed by induction on $\xi$.\bigskip

\noindent\bf The case $\xi\! =\! 1$\rm\bigskip

 We set $P_1\! :=\!\{\alpha\!\in\!\omega^\omega\mid\forall n\!\in\!\omega\ \ \alpha(2n)\! =\! 0\}$ and 
$\tau_1\! :=\!\tau$. The properties (1)$_1$-(3)$_1$ are clearly satisfied.\bigskip

\noindent (4)$_1$ Note that $(P_1,\tau_1)$ is homeomorphic to $(\omega^\omega ,\tau )$. As any nonempty open subset of $(\omega^\omega ,\tau )$ is homeomorphic to $(\omega^\omega ,\tau )$, 
$(U,\tau_1)$ is homeomorphic to $(\omega^\omega ,\tau )$. This gives $\varphi_{\xi ,U}$, which is nice since $\omega^\omega$ is closed in itself. This shows that we can take $G\! :=\! U$.\bigskip

\noindent (5)$_1$ As in (4)$_1$ we see that $(V,\tau_1)$ is homeomorphic to $(\omega^\omega ,\tau )$, and we can take $H\! :=\! V$.\bigskip

\noindent (6)$_1$ Note that $U$ is the disjoint union of a sequence $(C_n)_{n\in\omega}$ of nonempty clopen subsets of $(P_1,\tau_1)$. Let $(U_{1,n})_{n\in\omega}$ be a partition of 
$W\!\setminus\! P_1$ into clopen subsets of $(\omega^\omega ,\tau_1)$. As any nonempty open subset of $(P_1,\tau_1)$ or $(\omega^\omega ,\tau_1)$ is homeomorphic to 
$(\omega^\omega ,\tau )$, we can find homeomorphisms 
$$\varphi_0\! :\! (C_0,\tau_1)\!\rightarrow\! (\bigcup_{n>0}~U_{1,n},\tau_1)$$ 
and $\varphi_1\! :\! (\bigcup_{n>0}~C_n,\tau_1)\!\rightarrow\! (U_{1,0},\tau_1)$. As $C_0$ and 
$U_{1,0}$ are $\tau$-closed, $\varphi_0$ and $\varphi_1$ are nice. This shows that the gluing of $\varphi_0$ and $\varphi_1$ is a nice homeomorphism from $(U,\tau_1)$ onto 
$(W\!\setminus\! P_1,\tau_1)$. Thus we can take $G\! :=\! U$ and $K\! :=\! W\!\setminus\! P_1$.\bigskip

\noindent (7)$_1$ As in (6)$_1$ we write $V$ as the disjoint union of a sequence $(D_n)_{n\in\omega}$ of nonempty clopen subsets of $(\omega^\omega ,\tau_1)$. As the $(D_n,\tau_1)$'s are homeomorphic to $(\omega^\omega ,\tau_1)$, we can take $H\! :=\! V\!\setminus\! P_1$ and 
$L\! :=\! W\!\setminus\! P_1$.

\vfill\eject

\noindent\bf The induction step\rm\bigskip

 We assume that $1\! <\! \xi\! <\!\omega_1$ and that the assertion holds for each ordinal 
$\theta\! <\!\xi$. We fix a sequence of ordinals $(\xi_n)_{n\in\omega}$ containing each ordinal in 
$\xi\!\setminus\!\{ 0\}$ infinitely many times. We set
$$\begin{array}{ll}
& P_{\xi}\! =\!\Baire\!\times\! (\Pi_{i\in\omega}~\neg P_{\xi_i})\mbox{,}\cr
& \tau_{\xi}^{<}\! =\!\tau\!\times\! (\Pi_{i\in\omega}~\tau_{\xi_i})\mbox{,}\cr
& U_{\xi,n}\! =\!\Baire\!\times\! (\Pi_{i<n}~\neg P_{\xi_i})\!\times\! P_{\xi_n}\!\times\! 
(\Baire)^\omega ~~~~(n\!\in\!\omega ).
\end{array}$$
The family $\{U_{\xi,n}\mid n\!\in\!\omega\}$ is disjoint. We set $\sigma_0\! =\!\tau_\xi^<$ and 
$\sigma_{n+1}\! =\!\langle\sigma_n\cup\{ U_{\xi,n}\}\rangle$. It is easy to check that 
$U_{\xi,n}\!\in\!\bormone (\sigma_n)$. Applying Lemma 2.4 we get a topology $\tau_\xi\! :=\!\sigma_\infty$ such that\bigskip

(a) $\tau_\xi$ is zero-dimensional Polish,\smallskip

(b) ${\tau_\xi}_{\vert P_\xi}\! =\! {\tau_\xi^<}_{\vert P_\xi}$,\smallskip

\noindent and, for every $n\!\in\!\omega$,\smallskip

(c) ${\tau_\xi}_{\vert U_{\xi,n}}\! =\! {\tau_\xi^<}_{\vert U_{\xi,n}}$,\smallskip

(d) $U_{\xi,n}$ is $\tau_\xi$-clopen.\bigskip

\noindent We defined the topology $\tau_\xi$ on $(\Baire )^\omega$ instead of $\Baire$. However, since the spaces $\big( (\Baire )^\omega ,\tau^\omega )$ and $(\Baire , \tau )$ are homeomorphic we can replace the latter space by the former one in the proof. Since there is no danger of confusion we will write $\tau$ instead of $\tau^\omega$ to simplify the notation.\bigskip

\noindent (1)$_\xi$ Clearly, $\tau\!\subseteq\!\tau_{\xi}$. Note that $U_{\xi,n}\!\in\!\boraxi (\tau)$ for every $n\!\in\!\omega$ and $\tau_\xi^<\!\subseteq\!\boraxi (\tau )$, so that 
$\tau_\xi\!\subseteq\!\boraxi (\tau)$. Moreover, $(\Baire ,\tau_\xi )$ is clearly perfect.\bigskip

\noindent (2)$_\xi$ As $U_{\xi,n}$ is $\tau_\xi$-clopen, $P_\xi$ is $\tau_\xi$-closed. Note that 
${\tau_\xi}_{\vert P_\xi}\! =\! {\tau_\xi^<}_{\vert P_\xi}$ and $P_\xi$ contains no nonempty basic 
$\tau_\xi^<$-open set. This implies that $P_\xi$ is $\tau_\xi$-nowhere dense.\bigskip

\noindent (3)$_\xi$ Let $S\!\in\!\boraxi (\tau )$ be $\tau_\xi$-nonmeager in $P_\xi$. We may assume that $S\!\in\!\bormth (\tau)$ for some $\theta < \xi$. As 
${\tau_\xi}_{\vert P_\xi}\! =\! {\tau_\xi^<}_{\vert P_\xi}$ and $S$ has the Baire property with respect to the topology $\tau_\xi^<$ there exists a $\tau_\xi^<$-open set $V$ such that $S$ is $\tau_\xi^<$-comeager in $P_\xi\cap V$. Moreover, we may assume that $V$ has the following form:
$$V\! =\!\tilde V\!\times\! (\Pi_{i\leq k}~V_i)\!\times\! (\Baire)^\omega\mbox{,}$$
where $\tilde V\!\in\!\tau$, $V_i\!\in\!\tau_{\xi_i}$ and $V_i\!\subseteq\!\neg P_{\xi_i}$ for each  
$i\!\leq\! k$. The set $V^*\! =\!\tilde V\!\times\! (\Pi_{i\leq k}~V_i)\!\times\! (\Pi_{i>k}~\neg P_{\xi_i})$ 
is $\tau_\xi^<$-comeager in $V$ since $\neg P_{\xi_i}$ is $\tau_{\xi_i}$-comeager in $\Baire$ for every $i\!\in\!\omega$. As $P_\xi\cap V\! =\! V^*$, $S$ is $\tau_\xi^<$-comeager in $V^*$. Let 
$p\!\in\!\omega$ be such that $p\! >\! k$ and $\xi_p\!\geq\!\theta$. Define
$$\begin{array}{ll}
& \tau^*\! =\!\tau\!\times\! (\Pi_{i\not= p}~\tau_{\xi_i})\mbox{,}\cr
& Z\! =\!\tilde V\!\times\! V_0\!\times\!\cdots\!\times\! V_k\!\times\!\neg P_{\xi_{k+1}}\!\times\!\cdots\!\times\! \neg P_{\xi_{p-1}}\!\times\! (\Baire)^{\omega}\mbox{,}\cr
& \tau^{\sharp}\! =\!\tau\!\times\! (\Pi_{i<p}~\tau_{\xi_i})\!\times\!\tau\!\times\! 
(\Pi_{i>p}~\tau_{\xi_i}).
\end{array}$$
For $\alpha\!\in\!\Baire$ define a set $(\neg S)_\alpha$ by
$$(\neg S)_\alpha\! :=\!\{ (\tilde y, y_0, y_1,\dots , y_{p-1},y_{p+1},\dots)\!\in\!\Baire\mid 
(\tilde y, y_0, y_1,\dots, y_{p-1},\alpha ,y_{p+1},\dots )\!\in\!\neg S\}.$$

 Denote 
$S^*\! :=\!\{\alpha\!\in\!\Baire\mid (\neg S)_\alpha\mbox{ is }\tau^*\mbox{-nonmeager in }Z\}$. Note that $\neg S\!\in\!\borath (\tau )\!\subseteq\!\borath (\tau^{\sharp})$. By the Montgomery theorem (see 22.D in [K]), $S^*\!\in\!\borath (\tau )\!\subseteq\! {\bf\Sigma}_{\xi_p}^0(\tau )$. By the Kuratowski-Ulam theorem,  $S^*$ is $\tau_{\xi_p}$-meager in $\neg P_{\xi_p}$. Using the induction hypothesis, Condition (3)$_{\xi_p}$ implies that $S^*$ is $\tau_{\xi_p}$-meager in 
$P_{\xi_p}$. Using the Kuratowski-Ulam theorem again, we see that $S$ is $\tau_{\xi}^<$-comeager in the $\tau_\xi$-open set
$$W\! =\!\tilde V\!\times\! V_0\!\times\!\cdots\!\times\! V_k\!\times\!\neg P_{\xi_{k+1}}\!\times\!\cdots\! \times\!\neg P_{\xi_{p-1}}\!\times\! P_{\xi_p}\!\times\! (\Baire)^{\omega}.$$
As $W\!\subseteq\! U_{\xi,p}$, ${\tau_\xi}_{\vert W}\! =\! {\tau_\xi^<}_{\vert W}$ by (c), and consequently $S$ is $\tau_\xi$-comeager in $W$. Thus $S$ is $\tau_\xi$-nonmeager in 
$(\Baire)^\omega$ since $W$ is $\tau_\xi$-open.\bigskip

\noindent (4)$_\xi$ We first construct a $\tau_\xi$-dense open subset of $U$, which is the disjoint union of sets of the form
$$U^n\! :=\!\big( W^n\!\times\! (\Pi_{i<k_n}~W^n_i)\!\times\! (\Baire )^\omega\big)\cap P_\xi\! =\!
W^n\!\times\! (\Pi_{i<k_n}~W^n_i\!\setminus\! P_{\xi_i})\!\times\! (\Pi_{i\geq k_n}~\neg P_{\xi_i})\mbox{,}$$
where $W^n$ is a nonempty $\tau$-clopen set and $W^n_i$ is a nonempty $\tau_{\xi_i}$-clopen set. In order to do this, we fix an injective $\tau_\xi$-dense sequence $(x_n)_{n\in\omega}$ of $U$, which is possible since $(P_\xi ,\tau_\xi )$ is nonempty and perfect. We first choose $W^0$ and the $W^0_i$'s in such a way that $U^0$ is a proper $\tau_\xi$-clopen neighborhood of $x_0$ in $U$, which is possible since ${\tau_\xi}_{\vert P_\xi}\! =\! {\tau_\xi^<}_{\vert P_\xi}$. For the induction step, we choose $p_n$ minimal such that $x_{p_n}\!\notin\!\bigcup_{q\leq n}~U^q$. Then we choose $W^{n+1}$ and the 
$W^{n+1}_i$'s in such a way that $U^{n+1}$ is a proper $\tau_\xi$-clopen neighborhood of 
$x_{p_n}$ in $U\!\setminus\! (\bigcup_{q\leq n}~U^q)$.\bigskip

 There is a nice homeomorphism $\psi_n$ from $W^n$ onto $N_n\! :=\!\{\alpha\!\in\!\Baire\mid\alpha (0)\! =\! n\}$. The induction assumption gives,\smallskip
 
- for $i\! <\! k_n$, a $\tau_{\xi_i}$-dense $G_\delta$ subset $G^n_i$ of 
$W^n_i\!\setminus\! P_{\xi_i}$, and a nice $(\tau_{\xi_i},\tau )$-homeomorphism $\psi_{\xi_i,G^n_i}$ of $G^n_i$ onto $\omega^\omega$,\smallskip
 
- for $i\!\geq\! k_n$, a $\tau_{\xi_i}$-dense $G_\delta$ subset $G^n_i$ of $\neg P_{\xi_i}$, and a nice $(\tau_{\xi_i},\tau )$-homeomorphism $\psi_{\xi_i,G^n_i}$ of $G^n_i$ onto 
$\omega^\omega$.\smallskip

 By Lemma 2.1, the map $\psi_n\!\times\! (\Pi_{i\in\omega}~\psi_{\xi_i,G^n_i})$ is a nice 
$(\tau_\xi^<,\tau )$-homeomorphism from 
$$W^n\!\times\! (\Pi_{i\in\omega}~G^n_i)$$ 
onto $N_n\!\times\! (\Baire )^\omega$. If we set 
$G\! :=\!\bigcup_{n\in\omega}~\big( W^n\!\times\! (\Pi_{i\in\omega}~G^n_i)\big)$, then we get a nice 
$(\tau_\xi^<,\tau )$-homeomorphism from $G$ onto $\Baire$. We are done since 
${\tau_\xi}_{\vert P_\xi}\! =\! {\tau_\xi^<}_{\vert P_\xi}$.\bigskip

\noindent (5)$_\xi$ We essentially argue as in (4)$_\xi$. As $P_\xi$ is $\tau_\xi$-closed nowhere dense, we may assume that 
$$V\!\subseteq\!\neg P_\xi\! =\!\bigcup_{n\in\omega}~U_{\xi ,n}.$$ 
We first construct a $\tau_\xi$-dense open subset of $V\cap U_{\xi ,n}$, which is the disjoint union of sets of the form $V^{n,p}\! :=\! W^{n,p}\!\times\! (\Pi_{i<n}~W^{n,p}_i\!\setminus\! P_{\xi_i})\!\times\! 
(W^{n,p}_n\cap P_{\xi_n})\!\times\! (\Pi_{n<i<k_n^p}~W^{n,p}_i)\!\times\! (\Baire )^\omega$, where $W^{n,p}$ is a nonempty $\tau$-clopen set and $W^{n,p}_i$ is a nonempty $\tau_{\xi_i}$-clopen set. This is possible since ${\tau_\xi}_{\vert U_{\xi ,n}}\! =\! {\tau_\xi^<}_{\vert U_{\xi ,n}}$. We are done since $U_{\xi ,n}$ is $\tau_\xi$-clopen.\bigskip

\noindent (6)$_\xi$ As in (4)$_\xi$ we construct a $\tau_\xi$-dense open subset of $U$, which is the disjoint union of sets of the form 
$U^n\! :=\!\big( W^n\!\times\! (\Pi_{i<k_n}~W^n_i)\!\times\! (\Baire )^\omega\big)\cap P_\xi\! =\!
W^n\!\times\! (\Pi_{i<k_n}~W^n_i\!\setminus\! P_{\xi_i})\!\times\! (\Pi_{i\geq k_n}~\neg P_{\xi_i})$, 
where $W^n$ is a nonempty $\tau$-clopen set and $W^n_i$ is a nonempty $\tau_{\xi_i}$-clopen set. Recall also that 
$$U_{\xi,n}\! =\!\Baire\!\times\! (\Pi_{i<n}~\neg P_{\xi_i})\!\times\! P_{\xi_n}\!\times\! (\Baire)^\omega .$$ 
We also construct a $\tau_\xi$-dense open subset of $W$, which is the disjoint union of sets of the form
$$\pi^n\! :=\! Z^n\!\times\! (\Pi_{i<l_n}~Z^n_i\!\setminus\! P_{\xi_i})\!\times\! 
(Z^n_{l_n}\cap P_{\xi_{l_n}})\!\times\! (\Pi_{l_n<i<m_n}~Z^n_i)\!\times\! 
(\Baire )^\omega\!\subseteq\! U_{\xi ,l_n}\mbox{,}$$
where $Z^n$ is a nonempty $\tau$-clopen set and $Z^n_i$ is a nonempty $\tau_{\xi_i}$-clopen set. Let $(W^{0,p})_{p\in\omega}$ (respectively, 
$(Z^{0,p})_{p\in\omega}$) be a partition of $W^0$ (respectively, $Z^0$) into nonempty $\tau$-clopen sets. Using the facts that 
${\tau_\xi}_{\vert P_\xi}\! =\! {\tau_\xi^<}_{\vert P_\xi}$ and ${\tau_\xi}_{\vert U_{\xi ,n}}\! =\! {\tau_\xi^<}_{\vert U_{\xi ,n}}$, we will build\bigskip

- a nice $(\tau_\xi ,\tau_\xi )$-homeomorphism from a dense $G_\delta$ subset $G^{0,p}$ of 
$$U^{0,p}\! :=\! W^{0,p}\!\times\! (\Pi_{i<k_0}~W^0_i\!\setminus\! P_{\xi_i})\!\times\! 
(\Pi_{i\geq k_0}~\neg P_{\xi_i})$$ 
onto a dense $G_\delta$ subset $K^{0,p}$ of $\pi^{p+1}$. Then, using the fact that the 
$W^{0,p}$'s are $\tau$-clopen, the gluing of these homeomorphisms will be a nice 
$(\tau_\xi ,\tau_\xi )$-homeomorphism $\varphi_0$ from 
$G^0\! :=\!\bigcup_{p\in\omega}~G^{0,p}\!\subseteq\! U^0$ onto 
$K^0\! :=\!\bigcup_{p\in\omega}~K^{0,p}\!\subseteq\!\bigcup_{p>0}~\pi^p$.\smallskip

- a nice homeomorphism from a dense $G_\delta$ subset $G^{1,p}$ of $U^{p+1}$ onto a dense 
$G_\delta$ subset $K^{1,p}$ of $Z^{0,p}\!\times\! (\Pi_{i<l_0}~Z^0_i\!\setminus\! P_{\xi_i})\!\times\! 
(Z^0_{l_0}\cap P_{\xi_{l_0}})\!\times\! (\Pi_{l_0<i<m_0}~Z^0_i)\!\times\! (\Baire )^\omega$. Then the gluing of these homeomorphisms will be a nice $(\tau_\xi ,\tau_\xi )$-homeomorphism $\varphi_1$ from $G^1\! :=\!\bigcup_{p\in\omega}~G^{1,p}\!\subseteq\!\bigcup_{p>0}~U^p$ onto 
$K^1\! :=\!\bigcup_{p\in\omega}~K^{1,p}\!\subseteq\!\pi^0$.\bigskip

 The gluing of these two homeomorphisms will be a nice homeomorphism from 
$G\! :=\! G^0\cup G^1$ onto $K\! :=\! K^0\cup K^1$. The set $G^{0,p}$ (respectively, $K^{0,p}$) will be of the form $W^{0,p}\!\times\! (\Pi_{i\in\omega}~G^p_i)$ (respectively, $Z^{p+1}\!\times\! (\Pi_{i\in\omega}~K^p_i)$). Note first that there is a homeomorphism $\psi_p$ from $W^{0,p}$ onto $Z^{p+1}$. Then we build a permutation $i\!\mapsto\! j_i$ of the coordinates (with inverse $q\!\mapsto\! J_q$). This permutation is constructed in such a way that $\xi_{j_i}\! =\!\xi_i$, which will be possible since $(\xi_n)_{n\in\omega}$ contains each ordinal in $\xi\!\setminus\!\{ 0\}$ infinitely many times. If 
$i\! <\! m_{p+1}$ (respectively, $q\! <\! k_0$), then we choose $j_i\!\geq\! k_0$ (respectively, $J_q\!\geq\! m_{p+1}$), ensuring injectivity. For a remaining coordinate $q\!\notin\!\{ 0,...,k_0\! -\! 1\}\cup\{ j_l\mid l\! <\! m_{p+1}\}$, we choose $J_q\!\notin\!\{ 0,...,m_{p+1}\! -\! 1\}\cup\{ J_l\mid l\! <\! k_0\}$, ensuring that the map $q\!\mapsto\! J_q$ is a bijection from $\neg (\{ 0,...,k_0\! -\! 1\}\cup\{ j_l\mid l\! <\! m_{p+1}\} )$ onto 
$\neg\big(\{ 0,...,m_{p+1}\! -\! 1\}\cup\{ J_l\mid l\! <\! k_0\}\big)$. Then, using the induction assumption, we build our homeomorphism coordinate by coordinate, which means that $G^p_{j_i}$ will be homeomorphic to $K^p_i$. The induction assumption gives\bigskip
 
- for $i\! <\! l_{p+1}$, a $\tau_{\xi_{j_i}}$-dense $G_\delta$ subset $G^p_{j_i}$ of $\neg P_{\xi_{j_i}}$, a $\tau_{\xi_i}$-dense $G_\delta$ subset $K^p_i$ of 
$Z^{p+1}_i\!\setminus\! P_{\xi_i}$, and a nice $(\tau_{\xi_i},\tau_{\xi_i})$-homeomorphism 
$\psi_{\xi_i,G^p_{j_i},K^p_i}$ from $G^p_{j_i}$ onto $K^p_i$.\smallskip

- a $\tau_{\xi_{j_{l_{p+1}}}}$-dense $G_\delta$ subset $G^p_{j_{l_{p+1}}}$ of $\neg P_{\xi_{j_{l_{p+1}}}}$, a $\tau_{\xi_{l_{p+1}}}$-dense $G_\delta$ subset $K^p_{l_{p+1}}$ of $P_{\xi_{l_{p+1}}}$, and a nice $(\tau_{\xi_{l_{p+1}}},\tau_{\xi_{l_{p+1}}})$-homeomorphism 
$\varphi^{-1}_{\xi_{l_{p+1}},K^p_{l_{p+1}},G^p_{j_{l_{p+1}}}}$ from $G^p_{j_{l_{p+1}}}$ onto $K^p_{l_{p+1}}$.\smallskip

- for $l_{p+1}\! <\! i\! <\! m_{p+1}$, a $\tau_{\xi_{j_i}}$-dense $G_\delta$ subset 
$G^p_{j_i}$ of $\neg P_{\xi_{j_i}}$, a $\tau_{\xi_i}$-dense $G_\delta$ subset $K^p_i$ of 
$Z^{p+1}_i\!\setminus\! P_{\xi_i}$, and a nice $(\tau_{\xi_i},\tau_{\xi_i})$-homeomorphism 
$\psi_{\xi_i,G^p_{j_i},K^p_i}$ from $G^p_{j_i}$ onto $K^p_i$.\smallskip

- for $q\! <\! k_0$, a $\tau_{\xi_q}$-dense $G_\delta$ subset $G^p_q$ of 
$W^0_q\!\setminus\! P_{\xi_q}$, a $\tau_{\xi_{J_q}}$-dense $G_\delta$ subset $K^p_{J_q}$ of 
$\neg P_{\xi_{J_q}}$, and a nice $(\tau_{\xi_q},\tau_{\xi_q})$-homeomorphism 
$\psi_{\xi_q,G^p_q,K^p_{J_q}}$ from $G^p_q$ onto $K^p_{J_q}$.\smallskip

- for a remaining coordinate $q\!\notin\!\{ 0,...,k_0\! -\! 1\}\cup\{ j_l\mid l\! <\! m_{p+1}\}$, a $\tau_{\xi_q}$-dense $G_\delta$ subset $G^p_q$ of 
$\neg P_{\xi_q}$, a $\tau_{\xi_{J_q}}$-dense $G_\delta$ subset $K^p_{J_q}$ of $\neg P_{\xi_{J_q}}$, and a nice $(\tau_{\xi_q},\tau_{\xi_q})$-homeomorphism $\psi_{\xi_q,G^p_q,K^p_{J_q}}$ from $G^p_q$ onto $K^p_{J_q}$.\bigskip

 By Lemma 2.1, the product $\varphi^0_p$ of $\psi_p$ with these nice homeomorphisms is a nice 
$(\tau_\xi^<,\tau_\xi^<)$-homeomorphism from 
$G^{0,p}\! :=\! W^{0,p}\!\times\! (\Pi_{i\in\omega}~G^p_i)$ onto 
$K^{0,p}\! :=\! Z^{p+1}\!\times\! (\Pi_{i\in\omega}~K^p_i)$, as well as a $(\tau_\xi ,\tau_\xi )$-homeomorphism since ${\tau_\xi}_{\vert P_\xi}\! =\! {\tau_\xi^<}_{\vert P_\xi}$ and 
${\tau_\xi}_{\vert U_{\xi ,l_{p+1}}}\! =\! {\tau_\xi^<}_{\vert U_{\xi ,l_{p+1}}}$. As $G^0$ is the sum of the $G^{0,p}$'s, $G$ is a $\tau_\xi$-dense $G_\delta$ subset of $U^0$. Similarly, $K^0$ is a 
$\tau_\xi$-dense $G_\delta$ subset of $\bigcup_{p>0}~\pi^p$. Moreover, the gluing 
$\varphi^0$ of the $\varphi^0_p$'s is a $(\tau_\xi ,\tau_\xi )$-homeomorphism from $G^0$ onto $K^0$.\bigskip

 The construction of $\varphi^1$ is similar.\bigskip

\noindent (7)$_\xi$ We argue as in (6)$_\xi$.\hfill{$\square$}

\begin{lem} Let $1 \leq \xi < \omega_1$. Then there are disjoint families ${\cal F}_\xi$, 
${\cal G}_\xi$ of subsets of $\Baire$ and a topology $T_\xi$ on $\Baire$ such that\smallskip

(a)$_\xi$ $T_\xi$ is zero-dimensional perfect Polish and 
$\tau\!\subseteq\! T_\xi\!\subseteq\!\boraxi (\tau )$,\smallskip

(b)$_{\xi}$ ${\cal F}_\xi$ is $T_\xi$-dense, i.e., for any nonempty $T_\xi$-open set $V$, there is 
$F\!\in\! {\cal F}_\xi$ with $F\!\subseteq\! V$,\smallskip

\noindent and, for every $F\!\in\! {\cal F}_\xi$,\smallskip

(c)$_{\xi}$ $F$ is nonempty, $T_\xi$-nowhere dense, and in $\bormtwo (T_\xi )$,\smallskip

(d)$_{\xi}$ if $S\!\in\!\boraxi (\tau )$ is $T_\xi$-nonmeager in $F$, then $S$ is $T_\xi$-nonmeager in $\Baire$,\smallskip

(e)$_{\xi}$ there is a nice $(T_\xi ,\tau )$-homeomorphism $\varphi_F$ from $F$ onto $\Baire$,\smallskip

(f)$_{\xi}$ for any nonempty $T_\xi$-open sets $V,V'$, there are disjoint $G,G'\!\in\! {\cal G}_\xi$ with $G\!\subseteq\! V$, $G'\!\subseteq\! V'$, and there is a nice $(T_\xi ,T_\xi )$-homeomorphism 
$\varphi_{G,G'}$ from $G$ onto $G'$,\smallskip

\noindent and, for every $G\!\in\! {\cal G}_\xi$,\smallskip

(g)$_{\xi}$ $G$ is nonempty, $T_\xi$-nowhere dense, and in $\bormtwo (T_\xi )$,\smallskip

(h)$_{\xi}$ if $S\!\in\!\boraxi (\tau )$ is $T_\xi$-nonmeager in $G$, then $S$ is $T_\xi$-nonmeager in 
$\Baire$.\end{lem}

\noindent\bf Proof.\rm\ Let $P_\xi$ and $\tau_\xi$ be as in Lemma 2.6. We set 
$T_\xi\! =\! (\tau_\xi )^\omega$. Let $(U_n)_{n\in\omega}$ be a basis for the topology $T_\xi$ made of nonempty sets. For each $n\!\in\!\omega$, there is a finite sequence $(V^n_i)_{i<k_n}$ of nonempty $\tau_\xi$-open sets such that 
$(\Pi_{i<k_n}~V^n_i)\!\times\! (\Baire )^\omega\!\subseteq\! U_n$. Moreover, the sequence 
$(k_n)_{n\in\omega}$ is chosen to be strictly increasing. Lemma 2.6 provides\bigskip

- for $i\! <\! k_n$, a $\tau_\xi$-dense $G_\delta$ subset $H^n_i$ of $V^n_i\!\setminus\! P_\xi$ and a nice $(\tau_\xi ,\tau )$-homeomorphism $\psi_{\xi ,H^n_i}$,\smallskip

- a $\tau_\xi$-dense $G_\delta$ subset $G^n_{k_n}$ of $P_\xi$ and a nice 
$(\tau_\xi ,\tau )$-homeomorphism $\varphi_{\xi ,G^n_{k_n}}$,\smallskip

- for $i\! >\! k_n$, a $\tau_\xi$-dense $G_\delta$ subset $H^n_i$ of $\Baire$ and a nice 
$(\tau_\xi ,\tau )$-homeomorphism $\psi_{\xi ,H^n_i}$.

\vfill\eject

 We then put $F_n\! :=\! (\Pi_{i<k_n}~H^n_i)\!\times\! G^n_{k_n}\!\times\! (\Pi_{i>k_n}~H^n_i)$, so that $F_n\!\subseteq\! U_n$. We set ${\cal F}_\xi\! =\!\{ F_n\mid n\!\in\!\omega\}$. Then ${\cal F}_\xi$ is clearly a disjoint family and the properties (a)$_\xi$ and (b)$_\xi$ are obviously satisfied.\bigskip

\noindent (c)$_\xi$ As $P_\xi$ is $\tau_\xi$-nowhere dense, each $F_n$ is $T_\xi$-nowhere dense. Each $F_n$ is obviously also in $\bormtwo (T_\xi )$.\bigskip

\noindent (d)$_\xi$ Let $n\!\in\!\omega$ and $S\!\in\!\boraxi (\tau )$ be $T_\xi$-nonmeager in 
$F_n$. We define 
$$\begin{array}{ll}
& Z\! =\!\Pi_{i\not= k_n}~H^n_i\mbox{,}\cr
& T_\xi^*\! =\!\Pi_{i\not= k_n}~ {\tau_\xi}_{\vert H^n_i}\mbox{,}\cr
& \tilde T_\xi\! =\! (\Pi_{i<k_n}~ {\tau_\xi}_{\vert H^n_i})\!\times\!\tau\!\times\! 
(\Pi_{i>k_n}~ {\tau_\xi}_{\vert H^n_i}). 
\end{array}$$
If $\alpha\!\in\!\Baire$, then we denote 
$$S_\alpha\! :=\!\{ (y_0,\dots ,y_{k_n-1},y_{k_n+1},\dots )\!\in\!\Baire\mid 
(y_0,\dots ,y_{k_n-1},\alpha ,y_{k_n+1},\dots )\!\in\! S\} .$$ 
We set $S^*\! =\!\{\alpha\!\in\!\Baire\mid S_\alpha\mbox{ is }T_\xi^*\mbox{-nonmeager}\}$. By the Montgomery theorem,  $S^*\!\in\!\boraxi (\tau)$ since $S\!\in\!\boraxi (\tilde T_\xi )$. The set $S^*$ is $\tau_\xi$-nonmeager in $G^n_{k_n}$ by the Kuratowski-Ulam theorem, in $P_\xi$ also, and thus $S^*$ is $\tau_\xi$-nonmeager in $\Baire$. Using the Kuratowski-Ulam theorem again, we see that $S$ is $T_\xi$-nonmeager in 
$(\Pi_{i<k_n}~H^n_i)\!\times\!\Baire\!\times\! (\Pi_{i>k_n}~H^n_i)$, and thus in $(\Baire)^\omega$.\bigskip

\noindent (e)$_\xi$ We set $\varphi_F\! =\! (\Pi_{i<k_n}~\psi_{\xi ,H^n_i})\!\times\!
\varphi_{\xi ,G^n_{k_n}}\!\times\! (\Pi_{i>k_n}~\psi_{\xi ,H^n_i})$. The map $\varphi_F$ is clearly a 
$(T_\xi ,\tau )$-homeo-morphism from $F$ onto $(\Baire)^\omega$. It is nice by Lemma 2.1.\bigskip

 We now construct ${\cal G}_\xi$. For each $m\!\in\!\omega$, there are finite sequences 
$(V^m_i)_{i<k_m}$, $(W^m_i)_{i<l_m}$ of nonempty $\tau_\xi$-open sets such that 
$(\Pi_{i<k_m}~V^m_i)\!\times\! (\Baire )^\omega\!\subseteq\! U_{(m)_0}$ and 
${(\Pi_{i<l_m}~W^m_i)\!\times\! (\Baire )^\omega\!\subseteq\! U_{(m)_1}}$. Moreover, the sequences $(k_m)_{m\in\omega}$ and $(l_m)_{m\in\omega}$ are chosen to be strictly increasing and disjoint. Assume for example that $k_m\! <\! l_m$. Lemma 2.6 provides\bigskip

- for $i\! <\! k_m$, a $\tau_\xi$-dense $G_\delta$ subset $H^m_i$ of $V^m_i\!\setminus\! P_\xi$, a 
$\tau_\xi$-dense $G_\delta$ subset $L^m_i$ of $W^m_i\!\setminus\! P_\xi$, and a nice 
$(\tau_\xi ,\tau_\xi )$-homeomorphism $\psi_{\xi ,H^m_i,L^m_i}$,\smallskip

- a $\tau_\xi$-dense $G_\delta$ subset $G^m_{k_m}$ of $P_\xi$, a $\tau_\xi$-dense $G_\delta$ subset $K^m_{k_m}$ of $W^m_i\!\setminus\! P_\xi$, and a nice $(\tau_\xi ,\tau_\xi )$-homeomorphism $\varphi_{\xi ,G^m_{k_m},K^m_{k_m}}$,\smallskip

- for $k_m\! <\! i\! <\! l_m$, a $\tau_\xi$-dense $G_\delta$ subset $H^m_i$ of $\neg P_\xi$, a 
$\tau_\xi$-dense $G_\delta$ subset $L^m_i$ of $W^m_i\!\setminus\! P_\xi$, and a nice 
$(\tau_\xi ,\tau_\xi )$-homeomorphism $\psi_{\xi ,H^m_i,L^m_i}$,\smallskip

- a $\tau_\xi$-dense $G_\delta$ subset $K^m_{l_m}$ of $\neg P_\xi$, a $\tau_\xi$-dense 
$G_\delta$ subset $G^m_{l_m}$ of $P_\xi$, and a nice $(\tau_\xi ,\tau_\xi )$-homeomorphism 
$\varphi^{-1}_{\xi ,G^m_{l_m},K^m_{l_m}}$,\smallskip

- for $i\! >\! l_m$, a $\tau_\xi$-dense $G_\delta$ subset $H^m_i$ of $\neg P_\xi$, a $\tau_\xi$-dense $G_\delta$ subset $L^m_i$ of $\neg P_\xi$, and a nice $(\tau_\xi ,\tau_\xi )$-homeomorphism $\psi_{\xi ,H^m_i,L^m_i}$.\bigskip

\noindent We then put 
$$\begin{array}{ll}
& F'_m\! :=\! (\Pi_{i<k_m}~H^m_i)\!\times\! G^m_{k_m}\!\times\! (\Pi_{k_m<i<l_n}~H^m_i)
\!\times\! K^m_{l_m}\!\times\! (\Pi_{i>l_m}~H^m_i)\mbox{,}\cr
& G_m\! :=\! (\Pi_{i<k_m}~L^m_i)\!\times\! K^m_{k_m}\!\times\! (\Pi_{k_m<i<l_m}~L^m_i)
\!\times\! G^m_{l_m}\!\times\! (\Pi_{i>l_m}~L^m_i)\mbox{,}
\end{array}$$
so that $F'_m\!\times\! G_m\!\subseteq\! U_{(m)_0}\!\times\! U_{(m)_1}$. We set 
${\cal G}_\xi\! =\!\{ F'_m\mid m\!\in\!\omega\}\cup\{ G_m\mid m\!\in\!\omega\}$. Then ${\cal G}_\xi$ is clearly a disjoint family.

\vfill\eject

\noindent (f)$_\xi$ The map $\varphi_{F'_m,G_m}$ is by definition 
$$(\Pi_{i<k_m}~\psi_{\xi ,H^m_i,L^m_i})\!\times\!
\varphi_{\xi ,G^m_{k_m},K^m_{k_m}}\!\times\! (\Pi_{k_m<i<l_m}~\psi_{\xi ,H^m_i,L^m_i})\!\times\!
\varphi^{-1}_{\xi ,G^m_{l_m},K^m_{l_m}}\!\times\! (\Pi_{i>l_m}~\psi_{\xi ,H^m_i,L^m_i}).$$ 
Note that $\varphi_{F'_m,G_m}$ is clearly a $(T_\xi ,T_\xi )$-homeomorphism from $F'_m$ onto $G_m$. It is nice by Lemma 2.1.\bigskip

\noindent (g)$_\xi$ We argue as in (c)$_\xi$.\bigskip

\noindent (h)$_\xi$ We argue as in (d)$_\xi$.\hfill{$\square$}

\section{$\!\!\!\!\!\!$ Negative results}

\noindent\bf Proof of Theorem 1.1.\rm\ We apply Lemma 2.7 to the ordinal $\xi\! +\! 1$, which gives a  family $\mathcal F_{\xi+1}$ and a topology $T_{\xi+1}$ satisfying (a)$_{\xi+1}$-(e)$_{\xi+1}$. Let 
$(U_n\!\times\! V_n)_{n\in\omega}$ be a sequence of nonempty sets such that\bigskip 

- $U_n\!\in\! T_{\xi+1}$, $V_n$ is $\tau$-clopen,\smallskip

- $\{ U_n\!\times\! V_n\mid n\!\in\!\omega\}$ is a basis for the topology $T_{\xi+1}\!\times\!\tau$.\bigskip

\noindent For each $n\!\in\!\omega$ we find 
$F_n\!\in\! {\cal F}_{\xi+1}\!\setminus\!\{ F_q\mid q\! <\! n\}$ with $F_n\!\subseteq\! U_n$. By the property (e)$_{\xi+1}$ of $\mathcal F_{\xi+1}$ we find, for each $n\!\in\!\omega$, a nice 
$(T_{\xi+1},\tau )$-homeomorphism $f_n$ from $F_n$ onto $V_n$. We define 
$f\! :\!\bigcup_{n\in\omega}~F_n\!\rightarrow\!\Baire$ by $f(x)\! :=\! f_n(x)$ if $x\!\in\! F_n$. As 
$\mathcal F_{\xi+1}$ is a disjoint family, $f$ is well-defined. The graph of $f$ is 
$\boratwo (\tau\!\times\!\tau)$ since each $\mbox{Gr}(f_n)$ is $(\tau\!\times\!\tau)$-closed.\bigskip

 Suppose, towards a contradiction, that there exist, for $n\!\in\!\omega$, $C_n\! \in\!\boraxi (\tau )$ and  
$D_n\! \in\!\borel (\tau )$ such that $\neg\mbox{Gr}(f)\! =\!\bigcup_{n\in\omega}~C_n\!\times\! D_n$. By the Baire category theorem there is $n_0\!\in\!\omega$ such that $C_{n_0}$ is $T_{\xi+1}$-nonmeager and $D_{n_0}$ is $\tau$-nonmeager. As $C_{n_0}$ has the Baire property, we find a nonempty 
$T_{\xi+1}$-open set $O_1$ such that $C_{n_0}$ is $T_{\xi+1}$-comeager in $O_1$. Similarly, we find a $\tau$-open set $O_2$ such that $D_{n_0}$ is $\tau$-comeager in $O_2$.\bigskip

 Let $n\!\in\!\omega$ and $F_n\!\subseteq\! O_1$. Suppose that $C_{n_0}$ is not $T_{\xi+1}$-comeager in $F_n$. Then $O_1\!\setminus\! C_{n_0}$ is $T_{\xi+1}$-nonmeager in $F_n$. Note that 
$O_1\!\in\!\boraxp (\tau )$ and $C_{n_0}\!\in\!\boraxi (\tau )$. Therefore 
$O_1\!\setminus\! C_{n_0}\!\in\!\boraxp (\tau )$. Thus $O_1\!\setminus\! C_{n_0}$ is $T_{\xi+1}$-nonmeager in $\Baire$  by (d)$_{\xi+1}$. Consequently, $O_1\!\setminus\! C_{n_0}$ is 
$T_{\xi+1}$-nonmeager in $O_1$, a contradiction. Thus $C_{n_0}$ is $T_{\xi+1}$-comeager in $F_n$ for any $n\!\in\!\omega$ with $F_n\!\subseteq\! O_1$.\bigskip

 Find $n\!\in\!\omega$ such that $\mbox{Gr}(f_n)\!\subseteq\! O_1\!\times\! O_2$. Then $C_{n_0}$ is  $T_{\xi+1}$-comeager in $F_n$ and $D_{n_0}$ is $\tau$-comeager in $V_n$. As $f_n$ is a homeomorphism, $f_n^{-1}(V_n\cap D_{n_0})$ is $T_{\xi+1}$-comeager in $F_n$. As 
$F_n\!\in\!\bormtwo (T_{\xi+1})$ there exists 
$\alpha\!\in\!  f_n^{-1}(V_n\cap D_{n_0})\cap F_n\cap C_{n_0}$. This implies that 
$\big(\alpha ,f_n(\alpha )\big)\!\in\! C_{n_0}\!\times\! D_{n_0}$, a contradiction.\hfill{$\square$}\bigskip
 
\noindent\bf Proof of Theorem 1.2.\rm\ Apply Lemma 2.7 to the ordinal $\xi+1$, which gives a family ${\cal G}_{\xi+1}$ and a topology $T_{\xi+1}$ satisfying (a)$_{\xi+1}$-(h)$_{\xi+1}$. Let 
${\cal U}\! =\!\{ U_n\mid n\!\in\!\omega\}$ be a basis for the space $(\Baire ,T_{\xi+1})$ made of nonempty sets. For each $n\!\in\!\omega$ we find $T_{\xi+1}$-open sets $V_n$, $W_n$ such that 
$$V_n\!\times\! W_n\!\subseteq\! B_{\tau\times\tau}\big(\Delta (\Baire ),2^{-n}\big)\cap 
(U_n \times U_n)\!\setminus\!\Delta(\Baire )$$ 
(we use the standard metric on $(\Baire ,\tau )$).

\vfill\eject

 By the properties (f)$_{\xi+1}$ and (g)$_{\xi+1}$ of ${\cal G}_{\xi+1}$ we find, for each 
$n\!\in\!\omega$, sets $F_n$ and $H_n$ from ${\cal G}_{\xi+1}$ such that
$$(*)~~~~F_n\!\subseteq\! V_n\!\setminus\! (\bigcup_{j<n}~F_j\cup H_j)\ \wedge\ 
H_n\!\subseteq\! W_n\!\setminus\!\big( F_n\cup (\bigcup_{j<n}~F_j\cup H_j)\big).$$
Moreover, there is a nice $(T_{\xi+1},T_{\xi+1})$-homeomorphism $f_n$ from $F_n$ onto $H_n$. We set $${\cal G}\! =\!\bigcup\{\mbox{Gr}(f_n)\mid n\!\in\!\omega\} .$$ 
Now we check the desired properties.\bigskip

 As $\tau\!\subseteq\! T_{\xi+1}$, $\overline{\cal G}^{\tau\times\tau}\! =\! {\cal G}\cup\Delta(\Baire )$, by construction. Thus ${\cal G}$ is a difference of two $(\tau\!\times\!\tau )$-closed sets. As each $f_n$ is a homeomorphism, the property $(*)$ implies that $f$ is a partial injection with disjoint domain and range. In order to see that $\cal G$ has no $\borxi$-measurable countable coloring, we proceed by contradiction. Suppose that there are $\cal G$-discrete sets 
 $C_n\!\in\!\borxi (\tau )$ (a set $C$ is $\cal G$-{\bf discrete} if $C^2\cap {\cal G}\! =\!\emptyset$), for $n\!\in\!\omega$, such that $\Delta(\Baire )\!\subseteq\!\bigcup_{n \in \omega}~C_n^2$. By the Baire theorem there exists $n_0\!\in\!\omega$ such that $C_{n_0}$ is $T_{\xi+1}$-nonmeager. As $C_{n_0}$ has the Baire property, we find a nonempty $T_{\xi+1}$-open set $O$ such that 
$C_{n_0} \cap O$ is $T_{\xi+1}$-comeager in $O$.\bigskip

 Let $F\!\in\! {\cal G}_{\xi+1}$ with $F\!\subseteq\! O$. Suppose that $C_{n_0}$ is not 
$T_{\xi+1}$-comeager in $F$. Then $O\!\setminus\! C_{n_0}$ is $T_{\xi+1}$-nonmeager in $F$. Note that $O\!\in\!\boraxp (\tau )$ and $C_{n_0}\!\in\!\borxi (\tau )$. Therefore 
$O\!\setminus\! C_{n_0}\!\in\!\boraxp (\tau )$. Thus $O\!\setminus\! C_{n_0}$ is $T_{\xi+1}$-nonmeager in $\Baire$ by (h)$_{\xi+1}$. Consequently, $O\!\setminus\! C_{n_0}$ is $T_{\xi+1}$-nonmeager in $O$, a contradiction. Thus $C_{n_0}$ is $T_{\xi+1}$-comeager in $F$ for any 
$F\!\in\! {\cal G}_{\xi+1}$ with $F\!\subseteq\! O$.\bigskip

 Find $n\!\in\!\omega$ such that $\mbox{Gr}(f_n)\!\subseteq\! O^2$. Then $C_{n_0}$ is 
$T_{\xi+1}$-comeager in $F_n$ and in $H_n$. As $f_n$ is a homeomorphism, 
$f_n^{-1}(H_n\cap C_{n_0})$ is $T_{\xi+1}$-comeager in $F_n\!\in\!\bormtwo (T_{\xi+1})$. Thus there exists
$$\alpha\!\in\!  f_n^{-1}(H_n\cap C_{n_0})\cap F_n\cap C_{n_0}.$$ 
This implies that $\big(\alpha ,f_n(\alpha )\big)\!\in\! C_{n_0}^2$, a contradiction.\hfill{$\square$}

\section{$\!\!\!\!\!\!$ Positive results}

\bf (A) $\borxi$-measurable countable colorings\rm\bigskip

 In [L-Z], the following conjecture is made.\bigskip
 
\noindent\bf Conjecture\it\ Let $1\!\leq\!\xi\! <\!\omega_1$. Then there are\smallskip

- a $0$-dimensional Polish space $\mathbb{X}_\xi$,\smallskip

- an analytic relation $\mathbb{A}_{\xi}$ on $\mathbb{X}_\xi$\smallskip

\noindent such that for any ($0$-dimensional if $\xi\! =\! 1$) Polish space $X$, and for any analytic relation $A$ on $X$, exactly one of the following holds:\smallskip  

(a) there is a $\borxi$-measurable countable coloring of $A$ (i.e., a $\borxi$-measurable map $c\! :\! X\!\rightarrow\!\omega$ such that 
$A\!\subseteq\! (c\!\times\! c)^{-1}(\not= )$),\smallskip  

(b) there is a continuous map $f\! :\!\mathbb{X}_\xi\!\rightarrow\! X$ such that $\mathbb{A}_{\xi}\!\subseteq\! (f\!\times\! f)^{-1}(A)$.\rm

\vfill\eject

 This would be a $\borxi$-measurable version of the $\mathbb{G}_0$-dichotomy in [K-S-T]. This conjecture is proved for $\xi\!\leq\! 3$ in [L-Z]. Our goals here are the following. We want to give\smallskip
 
- a reasonable candidate for $\mathbb{A}_{\xi}$ in the general case,\smallskip

- an example for $\xi\! =\! 3$ that is much simpler than the one in [L-Z].\bigskip

 We set ${\bf\Pi}^0_0\!:=\!\borone$. The following result is proved in [M\'a] (see Theorem 4 and Lemma 13.(i)). 

\begin{thm} (M\'atrai) Let $1\!\leq\!\xi\! <\!\omega_1$. There are a true $\bormxi$ subset $P_\xi$ of 
$2^\omega$, and a Polish topology $\tau_\xi$ on $2^\omega$ such that\smallskip

(1)$_\xi$ $\tau_\xi$ is finer than the usual topology $\tau'$ on $2^\omega$,\smallskip

(2)$_\xi$ $P_\xi$ is $\tau_\xi$-closed and $\tau_\xi$-nowhere dense,\smallskip

(3)$_\xi$ if $G$ is a basic $\tau_\xi$-open set meeting $P_\xi$, and 
$D\!\in\!\bormlxi (2^\omega ,\tau')$ is such that  $D\cap P_\xi\cap G$ is comeager in 
$(P_\xi\cap G,{\tau_\xi}_{ \vert P_\xi\cap G})$, then there is a $\tau_\xi$-open set $G'$ such that $P_\xi\cap G'\! =\! P_\xi\cap G$ and $D\cap G'$ is comeager in 
$(G',{\tau_\xi}_{\vert G'})$.\end{thm} 

\noindent\bf Notation.\rm\ In the sequel $1\!\leq\!\xi\! <\!\omega_1$. Fix, for each $\xi$, an increasing sequence $(\eta_n)_{n\in\omega}$ of elements of $\xi$ (different from $0$ if $\xi\!\geq\! 2$) such that $\mbox{sup}_{n\in\omega}~(\eta_n\! +\! 1)\! =\!\xi$.\bigskip

\noindent $\bullet$ Let $<.,.>:\!\omega^2\!\rightarrow\!\omega$ be a bijection, defined for example by 
$<\! n,p\! >:=\! (\Sigma_{k\leq n+p}~k)\! +\! p$, whose inverse bijection is $q\!\mapsto\!\big((q)_{0},(q)_{1}\big)$.\bigskip

\noindent $\bullet$ If $u\!\in\! 2^{\leq\omega}$ and $n\!\in\!\omega$, then we define $(u)_n\!\in\! 2^{\leq\omega}$ by 
$(u)_n(p)\! :=\! u(<n,p>)$ if $<n,p><\!\vert u\vert$.\bigskip

\noindent $\bullet$ Let $(t_n)_{n\in\omega}$ be a dense sequence in $\omega^{<\omega}$ with $\vert t_n\vert\! =\! n$. For example,  let 
$(p_n)_{n\in\omega}$ be the sequence of prime numbers, and $I\! :\!\omega^{<\omega}\!\rightarrow\!\omega$ defined by $I(\emptyset )\! :=\! 1$, and 
$I(s)\! :=\! p_0^{s(0)+1}...p_{\vert s\vert -1}^{s(\vert s\vert -1)+1}$ if $s\!\not=\!\emptyset$. Note that $I$ is one-to-one, so that there is an increasing bijection 
$i\! :\! I[\omega^{<\omega}]\!\rightarrow\!\omega$. Set $\psi\! :=\! (i\circ I)^{-1}\! :\!\omega\!\rightarrow\!\omega^{<\omega}$, so that $\psi$ is a bijection. Note that $|\psi (n)|\!\leq\! n$ if $n\!\in\!\omega$. Indeed, 
$$I\big(\psi (n)\vert 0\big)\! <\! I\big(\psi (n)\vert 1\big)\! <\! ...\! <\! I\big(\psi (n)\big)\mbox{,}$$ 
so that $(b\circ I)\big(\psi (n)\vert 0\big)\! <\! (b\circ I)\big(\psi (n)\vert 1\big)\! <\! ...\! <\! (b\circ I)\big(\psi (n)\big)\! =\! n$. As $|\psi (n)|\!\leq\! n$, we can define $t_{n}\! :=\!\psi (n)0^{n-|\psi (n)|}$, and $(t_n)_{n\in\omega}$ is suitable.\bigskip

\noindent $\bullet$ Theorem 4.1 gives $P_\xi$ and $\tau_\xi$. Let $Q_\xi\! :=\! 2\!\times\! P_\xi$, 
$T_\xi\! :=\!\mbox{discrete}\!\times\!\tau_\xi$, and $T_\xi^<\! :=\!\Pi_{i\in\omega}~T_{\eta_i}$ if $\xi\!\geq\! 2$.\bigskip

\noindent $\bullet$ $(W_{\xi ,n})_{n\in\omega}$ is a sequence of nonempty $T_\xi$-open sets.\bigskip

\noindent $\bullet$ $S_i\! :=\! Q_{\eta_i}\cup\bigcup_{n\in\omega}~W_{\eta_i,n}$ (for $i\!\in\!\omega$), and $S\! :=\!\Pi_{i\in\omega}~S_i$, so that 
$S\!\in\!\bormtwo (T_\xi^<)$ is a Polish space.\bigskip

\noindent $\bullet$ If $\xi\!\geq\! 2$, then we set\bigskip 

\leftline{$\mathbb{K}_\xi\! :=\!\bigcup_{n\in\omega}~
\Big\{ (\alpha ,\beta )\!\in\! 2^\omega\!\times\! 2^\omega\mid
\big(\forall i\! <\! n~~(\alpha )_i\! =\! (\beta )_i\!\in\! W_{\eta_i,t_n(i)}\big)~\wedge$}\smallskip

\rightline{$\big(\exists\gamma\!\in\! P_{\eta_n}~~
\big( (\alpha )_n,(\beta )_n\big)\! =\! (0\gamma ,1\gamma )\big)~\wedge ~
\big(\forall i\! >\! n~~(\alpha )_i\! =\! (\beta )_i\big)\Big\}\mbox{,}$}

\begin{lem} Let $2\!\leq\!\xi\!\leq\!\omega_1$. We assume that $Q_{\eta_i}\!\subseteq\!\overline{\bigcup_{n\in\omega}~W_{\eta_i,n}}^{T_{\eta_i}}$ for each 
$i\!\in\!\omega$. Then any $\mathbb{K}_\xi$-discrete $\boraxi$ subset $C$ of $(S,\tau')$ is $T_\xi^<$-meager in $S$.\end{lem} 

\noindent\bf Proof.\rm\ We may assume that $C$ is $\bormlxi$. We argue by contradiction. This gives $n\!\in\!\omega$ with $C\!\in\!\bormen$, a basic 
$T_\xi^<$-open set $O$ such that $C\cap O$ is $T_\xi^<$-comeager in $O\cap S\!\not=\!\emptyset$, $l\!\geq\! n$, and a sequence $(O_i)_{i<l}$ with 
$O_i\!\in\! T_{\eta_i}$ and ${O\! =\!\{\alpha\!\in\! 2^\omega\mid\forall i\! <\! l~~(\alpha )_i\!\in\! O_i\}}$. The assumption gives, for each $i\! <\! l$, $n_i\!\in\!\omega$ such that $O_i\cap W_{\eta_i,n_i}\!\not=\!\emptyset$. Let $m\!\geq\! l$ such that $t_m(i)\! =\! n_i$ for each $i\! <\! l$, and  
$$U\! :=\!\Big\{\alpha\!\in\! S\mid\forall i\! <\! l~~(\alpha )_i\!\in\! O_i\ \wedge\ 
\forall i\! <\! m~~(\alpha )_i\!\in\! W_{\eta_i,t_m(i)}\Big\}\mbox{,}$$ 
which is a nonempty $T_\xi^<$-open subset of $S$. In particular, $C\cap U$ is $T_\xi^<$-comeager in $U$. We set 
$$V\! :=\!\Big\{ (\alpha_i)_{i\not= m}\!\in\!\pi_{i\not= m}~S_i\mid\forall i\! <\! l~~
\alpha_i\!\in\! O_i\ \wedge\ \forall i\! <\! m~~\alpha_i\!\in\! W_{\eta_i,t_m(i)}\Big\}\mbox{,}$$ 
so that, up to a permutation of coordinates, $U\!\equiv\! S_m\!\times\! V$. We also set 
$$C'\! :=\!\Big\{\alpha\!\in\! S_m\mid\big( C\cap (S_m\!\times\! V)\big)_\alpha\mbox{ is }
\pi_{i\not= m}~T_{\eta_i}\mbox{-comeager in }V\Big\} .$$ 
By the Kuratowski-Ulam theorem, $C'$ is $T_{\eta_m}$-comeager in $S_m$ (see 8.41 in [K]). Write $C\! =\! D\cap S$, where $D\!\in\!\bormen (2^\omega )$. Note that $C'\! :=\! S_m\cap\Big\{\alpha\!\in\! 2^\omega\mid\big( D\cap (2^\omega\!\times\! V)\big)_\alpha\mbox{ is }\pi_{i\not= m}~T_{\eta_i}
\mbox{-comeager in }V\Big\}$. As $m\!\geq\! n$ and $\pi_{i\not= m}~T_{\eta_i}$ is finer than the usual topology, 
$D\cap (2^\omega\!\times\! V)\!\in\!\bormem (2^\omega ,\tau'\!\times\! (\pi_{i\not= m}~T_{\eta_i})_{\vert V})$. By the Montgomery theorem, $C'$ is 
$\bormem (S_m,\tau')$ (see 22.22 in [K]).\bigskip

 The set $C'$ cannot be $T_{\eta_m}$-comeager in $Q_{\eta_m}\cap N_0$ and $Q_{\eta_m}\cap N_1$. Indeed, we argue by contradiction to see that. We set $h_0(\alpha )\! :=<1\! -\!\alpha (0),\alpha (1),\alpha (2),...>$. As 
${h_0}_{\vert Q_{\eta_m}\cap N_0}$ is a $T_{\eta_m}$-homeomorphism, $C'\cap {h_0}_{\vert Q_{\eta_m}\cap N_0}^{-1}(C'\cap Q_{\eta_m}\cap N_1)$ is $T_{\eta_m}$-comeager in $Q_{\eta_m}\cap N_0$, and if $0\gamma$ is in it, then $1\gamma\!\in\! C'$, which gives 
$\delta\!\in\! (C\cap U)_{0\gamma}\cap (C\cap U)_{1\gamma}$ and contradicts the $\mathbb{K}_\xi$-discreteness of $C$.\bigskip

 Assume for example that $C'$ is not $T_{\eta_m}$-comeager in $Q_{\eta_m}\cap N_0$. Then 
$\neg C'$ is $T_{\eta_m}$-non meager in $Q_{\eta_m}$. As $C'$ is $\bormem (S_m,\tau')$, there is a sequence $(C_j)_{j\in\omega}$ of ${\bf\Pi}^0_{<\eta_m}(2^\omega )$ sets such that 
$$S_m\!\setminus\! C'\! =\!\bigcup_{j\in\omega}~C_j\cap S_m.$$ 
This gives $j\!\in\!\omega$ such that $C_j\cap Q_{\eta_m}$ is $T_{\eta_m}$-non meager in 
$Q_{\eta_m}$, and a basic $T_{\eta_m}$-open set $O$ such that $C_j\cap Q_{\eta_m}\cap O$ is $T_{\eta_m}$-comeager in 
$Q_{\eta_m}\cap O\!\not=\!\emptyset$.\bigskip

 The set $O$ is of the form $\{\varepsilon\}\!\times\! G$, where $\varepsilon\!\in\! 2$ and $G$ is a basic $\tau_{\eta_m}$-open set. Let $S\! :\! N_\varepsilon\!\rightarrow\! 2^\omega$ be the map defined by $S(\varepsilon\alpha )\! :=\!\alpha$. Note that $S$ is a $\tau'$-$\tau'$ and 
$T_\xi$-$\tau_\xi$ homeomorphism. In particular, 
$E\! :=\!\{\alpha\!\in\! 2^\omega\mid\varepsilon\alpha\!\in\! C_j\}$ is $\tau'$-${\bf\Pi}^0_{<\eta_m}$ and $E\cap P_{\eta_m}\cap G$ is comeager in 
$(P_{\eta_m}\cap G,{\tau_{\eta_m}}_{ \vert P_{\eta_m}\cap G})$. Theorem 4.1.(3) gives a 
$\tau_{\eta_m}$-open set $G'$ such that $P_{\eta_m}\cap G'\! =\! P_{\eta_m}\cap G$ and 
$E\cap G'$ is comeager in $(G',{\tau_{\eta_m}}_{\vert G'})$. Now $O'\! :=\!\{\varepsilon\}\!\times\! G'$ is a $T_{\eta_m}$-open set such that $Q_{\eta_m}\cap O'\! =\! Q_{\eta_m}\cap O$ and $C_j\cap O'$ is $T_{\eta_m}$-comeager in $O'$. The assumption gives $n\!\in\!\omega$ such that 
$W_{\eta_m,n}\cap O'\!\not=\!\emptyset$. Note that $C_j\cap W_{\eta_m,n}\cap O'$ is 
$T_{\eta_m}$-comeager in $W_{\eta_m,n}\cap O'$, so that $\neg C'$ is $T_{\eta_m}$-non meager in $S_m$, which is absurd.\hfill{$\square$}

\begin{cor} Let $2\!\leq\!\xi\!\leq\!\omega_1$. We assume that $Q_{\eta_i}\!\subseteq\!\overline{\bigcup_{n\in\omega}~W_{\eta_i,n}}^{T_{\eta_i}}$ for each 
$i\!\in\!\omega$. Then\smallskip

(a) there is no $\borxi$-measurable map $c\! :\! 2^\omega\!\rightarrow\!\omega$ such that 
$\mathbb{K}_\xi\!\subseteq\! (c\!\times\! c)^{-1}(\not= )$,\smallskip

(b) if $\mathbb{X}_\xi\!\in\!\bormtwo (2^\omega )$ and $\mathbb{K}_\xi\!\subseteq\!\mathbb{X}_\xi^2$, then there is no $\borxi$-measurable map 
$c\! :\!\mathbb{X}_\xi\!\rightarrow\!\omega$ such that $\mathbb{K}_\xi\!\subseteq\! (c\!\times\! c)^{-1}(\not= )$.\end{cor}

\noindent\bf Proof.\rm\ (a) We just have to apply Lemma 4.2.\bigskip

\noindent (b) We argue by contradiction. This gives a partition $(C_k)_{k\in\omega}$ of 
$\mathbb{X}_\xi$ into $\mathbb{K}_\xi$-discrete $\borxi (\mathbb{X}_\xi )$ sets. We set 
$D_0\! :=\! 2^\omega\!\setminus\!\mathbb{X}_\xi$, and choose $D_{k+1}\!\in\!\boraxi (2^\omega )$ such that $C_k\! =\! D_{k+1}\cap\mathbb{X}_\xi$. Then $(D_k)_{k\in\omega}$ is a covering of 
$2^\omega$ into $\mathbb{K}_\xi$-discrete $\boraxi$ sets. It remains to apply the reduction property of the class $\boraxi$ to contradict (a).\hfill{$\square$}\bigskip

\noindent\bf The case $\xi\! =\! 2$\rm\bigskip

\noindent\bf Example.\rm\ Let $\alpha\!\mapsto\!\alpha^*$ be the shift map on $2^\omega$: 
$\alpha^*(j)\! :=\!\alpha (j\! +\! 1)$. Then we set\bigskip 

\leftline{$\mathbb{A}_2\! :=\!\bigcup_{n\in\omega}~
\Big\{ (\alpha ,\beta )\!\in\! 2^\omega\!\times\! 2^\omega\mid
\big(\forall i\! <\! n~~(\alpha )_i\! =\! (\beta )_i\ \wedge\ 0^{t_n(i)}1\!\subseteq\! (\alpha )_i^*\big)~\wedge$}\smallskip

\rightline{$\big( (\alpha )_n,(\beta )_n\big)\! =\! (0^\infty ,10^\infty )~\wedge ~
\big(\forall i\! >\! n~~(\alpha )_i\! =\! (\beta )_i\big)\Big\} .$}

\begin{thm} $(2^\omega ,\mathbb{A}_2)$ satisfies the conjecture.\end{thm}

\noindent\bf Proof.\rm\ We set $P_1\! :=\!\{ 0^\infty\}$ and $\tau_1\! :=\!\tau'$, so that $P_1$ and 
$\tau_1$ satisfy the properties of Theorem 4.1. We also set 
$W_{1,n}\! :=\! N_{0^{n+1}1}\cup N_{10^n1}$, so that $(W_{1,n})_{n\in\omega}$ is a sequence of nonempty $T_1$-open sets satisfying the assumption of Corollary 4.3, so that 
$\mathbb{A}_2\! =\!\mathbb{K}_2$ satisfies its conclusions. In particular, (a) and (b) cannot hold simultaneously.\bigskip

 We define, for $(\varepsilon ,n)\!\in\! 2\!\times\!\omega$, $K_n^\varepsilon\! :=\!\{\alpha\!\in\! 2^\omega\mid\forall i\! <\! n~~0^{t_n(i)}1\!\subseteq\! (\alpha )_i^*\ \wedge\ (\alpha )_n(0)\! =\!\varepsilon\}$, and also $C_n^\varepsilon\! :=\! K_n^\varepsilon\!\setminus\!\big(\bigcup_{n<k}~K_k^0\cup K_k^1\big)$, so that $C_n^\varepsilon$ is closed, the $C^\varepsilon_n$'s are pairwise disjoint, and 
$\mathbb{A}_2\!\subseteq\!\bigcup_{n\in\omega}~C_n^0\!\times\! C_n^1$. We set, for each 
$p,q\!\in\!\omega$, 
$$O^p_q\! :=\!\left\{\!\!\!\!\!\!\!
\begin{array}{ll}
& K_n^\varepsilon\!\setminus\!\big(\bigcup_{n<k\leq q}~K_k^0\cup K_k^1\big)
\mbox{ if }p\! =\! 2n\! +\!\varepsilon\!\leq\! 2q\! +\! 1\mbox{,}\cr
& 2^\omega\!\setminus\! (\bigcup_{p'\leq 2q+1}~O^{p'}_q)\mbox{ if }p\! =\! 2q\! +\! 2\mbox{,}\cr
& \emptyset\mbox{ if }p\!\geq\! 2q\! +\! 3\mbox{,}
\end{array}
\right.$$
so that $(O^p_q)_{p\in\omega}$ is a covering of $2^\omega$ into clopen sets. Assume that 
$p\! =\! 2n\! +\!\varepsilon\!\not=\! p'\! =\! 2n'\! +\!\varepsilon'\!\leq\! 2q\! +\! 1$ and 
$\alpha\!\in\! O^p_q\cap O^{p'}_q$, so that $n,n'\!\leq\! q$. As 
$\alpha\!\in\! K^\varepsilon_n\cap K^{\varepsilon'}_{n'}$, $n\!\not=\! n'$ and for example $n\! <\! n'$, which is absurd. Thus $(O^p_q)_{p\in\omega}$ is a partition of $2^\omega$.\bigskip

\noindent (a) Assume that $q\! <\! n$. Note that $C^0_n\cup C^1_n$ is contained in or disjoint from each set of the form $K^\varepsilon_k$ with $k\!\leq\! q$. By disjonction, there is at most one couple 
$(\varepsilon ,r)$ such that $2r\! +\!\varepsilon\!\leq\! 2q\! +\! 1$ and 
$C^0_n\cup C^1_n\!\subseteq\! O^{2r+\varepsilon}_q$. If it does not exist, then 
$C^0_n\cup C^1_n\!\subseteq\! O^{2q+2}_q$.\bigskip

\noindent (b) Assume that $q\!\geq\! n$. Note that $C^\varepsilon_n\!\subseteq\! K^\varepsilon_n$. As 
$q\!\geq\! n$, $p\! :=\! 2n\! +\!\varepsilon\!\leq\! 2q\! +\! 1$. Thus $C^\varepsilon_n\!\subseteq\! O^p_q$.\bigskip

 It remains to apply Proposition 4.6 in [L-Z] to see that (a) or (b) holds.\hfill{$\square$}
 
\vfill\eject

\noindent\bf The case $\xi\! =\! 3$\rm\bigskip

\noindent\bf Example.\rm\  Let $(s_n)_{n\in\omega}$ be a dense sequence in $2^{<\omega}$ with 
$\vert s_n\vert\! =\! n$. For example, let $\phi\! :\!\omega\!\rightarrow\! 2^{<\omega}$ be a natural bijection. More specifically, $\phi (0)\! :=\!\emptyset$ is the sequence of length $0$, $\phi (1)\! :=\! 0$, 
$\phi (2)\! :=\! 1$ are the sequences of length $1$, and so on. Note that $|\phi (n)|\!\leq\! n$ if 
$n\!\in\!\omega$. Let $n\!\in\!\omega$. As $|\phi (n)|\!\leq\! n$, we can define 
$s_n\! :=\!\phi (n)0^{n-|\psi (n)|}$. We set $P_2\! :=\!\{\alpha\!\in\! 2^\omega\mid\forall p\!\in\!\omega ~~\exists q\!\geq\! p~~\alpha (q)\! =\! 1\}$, and\bigskip

\leftline{$\mathbb{A}_3\! :=\!\bigcup_{n\in\omega}~
\Big\{ (\alpha ,\beta )\!\in\! 2^\omega\!\times\! 2^\omega\mid
\big(\forall i\! <\! n~~(\alpha )_i\! =\! (\beta )_i\! =\! s_{t_n(i)}10^\infty\big)~\wedge$}\smallskip

\rightline{$\big(\exists\gamma\!\in\! P_2~~
\big( (\alpha )_n,(\beta )_n\big)\! =\! (0\gamma ,1\gamma )\big)~\wedge ~
\forall i\! >\! n~~(\alpha )_i\! =\! (\beta )_i\Big\} .$}\bigskip

\noindent We will see that $\mathbb{A}_3$, together with a suitable $\bormtwo$ subset 
$\mathbb{X}_3$ of $2^\omega$, satisfies the conjecture. The topology $\tau_2$ makes the countably many singletons of $\neg P_2$ open. Then $P_2$ is a true $\bormtwo$ subset of 
$2^\omega$ (see 23.A in [K]), $\tau_2$ is Polish finer than $\tau'$, $P_2$ is closed nowhere dense for $\tau_2$ since $\tau_2$ coincides with $\tau'$ on $P_2$ and $\neg P_2$ is $\tau'$-dense, and 4.1.(3) is satisfied since a basic $\tau_2$-open set meeting $P_2$ is a basic $\tau'$-clopen set and $P_2$ is $\tau'$-comeager. Thus $P_2$ and $\tau_2$ satisfy the properties of Theorem 4.1. We set 
$$W_{2,n}\! :=\!\{ s_n10^\infty\} .$$ 
Then $Q_2\!\subseteq\!\overline{\bigcup_{n\in\omega}~W_{2,n}}^{T_2}$ since $(s_n)_{n\in\omega}$ is dense. This shows that $\mathbb{A}_3\! =\!\mathbb{K}_3$ satisfies the conclusions of Corollary 4.3. In particular, (a) and (b) cannot hold simultaneously. In order to prove that (a) or (b) holds, we simply indicate the modifications to make to Section 5 in [L-Z]. We just need to prove the right lemmas since the final construction is the same.

\begin{lem} (a) Let $n\!\in\!\omega$ and $i\! <\! n$. Then $t_n(i)\! <\! n\! -\! i$.\smallskip

(b) The map $M\! :\!\{ s_{t_n(i)}10^\infty\mid n\!\in\!\omega\ \wedge\ i\! <\! n\}\!\rightarrow\!\omega$, defined by $M(\alpha )\! :=\!\mbox{max}\{ p\!\in\!\omega\mid\alpha (p)\! =\! 1\}$, is one-to-one.\end{lem}

\noindent\bf Proof.\rm\ (a) Recall the map $\psi$ defined after Theorem 4.1. It is enough to prove that 
$\psi (n)(i)\! <\! n\! -\! i$ if $i\! <\!\vert\psi (n)\vert$. We argue by induction on $n$, and the result is clear for $n\! =\! 0$. We may assume that $\psi (n)(i)\! =\! q\! +\! 1$ for some natural number $q$. We define 
$t\!\in\!\omega^{<\omega}$ by $t(i)\! :=\! q$, and $t(j)\! :=\!\psi (n)(j)$ if $j\!\not=\! i$. Let 
$p\!\in\!\omega$ with $\psi (p)\! =\! t$. Note that 
$I\big(\psi (p)\big)\! <\!I\big(\psi (n)\big)$, so that $p\! <\! n$. The induction assumption implies that 
$q\! =\!\psi (p)(i)\! <\! p\! -\! i$, so that $\psi (n)(i)\! =\! q\! +\! 1\!\leq\! p\! -\! i\! <\! n\! -\! i$.\bigskip 

\noindent (b) Assume that $M(\alpha )\! =\! M(\alpha')$. Let $n,n',i,i'$ with 
$\alpha\! =\! s_{t_n(i)}10^\infty$ and $\alpha'\! =\! s_{t_{n'}(i')}10^\infty$. Then  
$t_n(i)\! =\!\vert s_{t_n(i)}\vert\! =\! M(\alpha )\! =\! M(\alpha')\! =\! t_{n'}(i')$, so that $\alpha\! =\!\alpha'$.
\hfill{$\square$}\bigskip

\noindent\bf Notation.\rm\ If $\emptyset\!\not=\! u\!\in\! 2^{<\omega}$, then 
$u^m\! :=\! u\vert (\vert u\vert\! -\! 1)$.\bigskip

 The following notion is technical but crucial.

\begin{defi} We say that $u\!\in\! 2^{<\omega}$ is $placed$ if\smallskip

(a) $u\!\not=\!\emptyset$,\smallskip

(b) $\forall i\! <\! (\vert u^m\vert )_0~~(u)_i\!\subseteq\! s_{t_{(\vert u^m\vert )_0}(i)}10^\infty$,\smallskip

(c) $u(\vert u^m\vert )\! =\! 1$ if $(\vert u^m\vert )_1\! >\! 0$.\end{defi}

 We are now ready to define 
$$\mathbb{X}_3\! :=\!\big\{\alpha\!\in\! 2^\omega\mid\forall n\!\in\!\omega ~~
\exists p\!\geq\! n~~\alpha\vert p\mbox{ is placed}\big\} .$$ 
Note that $\mathbb{X}_3$ is a $\bormtwo$ subset of $2^\omega$. In particular, $\mathbb{X}_3$ is a 
$0$-dimensional Polish space.

\begin{lem} (a) The set $\mathbb{A}_3$ is a $\borathree$ (and thus analytic) relation on 
$\mathbb{X}_3$.\smallskip

(b) $(\mathbb{X}_3,\mathbb{A}_3)\not\preceq_{\borthree}\big(\omega ,\neg\Delta (\omega )\big)$.
\end{lem}

\noindent\bf Proof.\rm\ (a) $\mathbb{A}_3$ is clearly a $\borathree$ relation on $2^\omega$. So it is enough to see that it is a relation on $\mathbb{X}_3$. Fix $(\alpha ,\beta )\!\in\!\mathbb{A}_3$ (which defines a natural number $n$). Choose an infinite sequence $(p_k)_{k\in\omega}$ of natural numbers such that $(\alpha )_n(p_k)\! =\! (\beta )_n(p_k)\! =\! 1$. Then $\alpha\vert (<n,p_k>\! +1)$ and 
$\beta\vert (<n,p_k>\! +1)$ are placed, so that $\alpha ,\beta\!\in\!\mathbb{X}_3$.\bigskip

\noindent (b) This comes from Corollary 4.3.\hfill{$\square$}

\begin{lem}  Let $n\!\in\!\omega$, $\alpha\!\in\! 2^\omega$ such that $(\alpha )_i\! =\! s_{t_n(i)}10^\infty$ for each $i\! <\! n$, and $p\! > <n,0>$ such that $\alpha\vert p$ is placed. Then $(p\! -\! 1)_0\!\geq\! n$.\end{lem}

\noindent\bf Proof.\rm\ We argue by contradiction. As $p\! -\! 1\!\geq <n,0>$, 
$(p\! -\! 1)_0\! +\! (p\! -\! 1)_1\!\geq\! n\! +\! 0\! =\! n$. Thus $(p\! -\! 1)_1\!\geq\! n\! -\! (p\! -\! 1)_0\! >\! 0$. As 
$\alpha\vert p$ is placed, $\alpha (p\! -\! 1)\! =\! 1$. But 
$$\alpha (p\! -\! 1)\! =\!\alpha (<(p\! -\! 1)_0,(p\! -\! 1)_1>)\! =\! (\alpha )_{(p\! -\! 1)_0}\big( (p\! -\! 1)_1\big)
\! =\!\big( s_{t_n( (p\! -\! 1)_0)}10^\infty\big)\big( (p\! -\! 1)_1\big) .$$
By Lemma 4.5.(a), we get $(p\! -\! 1)_1\! <\! n\! -\! (p\! -\! 1)_0$, which is absurd.\hfill{$\square$}

\begin{defi} Let $u\!\in\! 2^{<\omega}$ and $l\!\in\!\omega$.\smallskip

(a) If $u$ is placed, then we will consider\smallskip

\noindent $\bullet$ the natural number $l(u)\! :=\! (\vert u^m\vert )_0$\smallskip

\noindent $\bullet$ the sequence $u^{l(u)}\!\in\! 2^{\vert u\vert}\!\setminus\!\{ u\}$ defined by 
$u^{l(u)}(m)\! :=\! u(m)$ exactly when $m\!\not= <l(u),0>$. Note that $u^{l(u)}$ is placed,  
$l(u^{l(u)})\! =\! l(u)$ and $(u^{l(u)})^{l(u)}\! =\! u$\smallskip

\noindent $\bullet$ the digit $\varepsilon (u)\! :=\! u(<l(u),0>)$. Note that 
$\varepsilon (u^{l(u)})\! =\! 1\! -\!\varepsilon (u)$.\smallskip

(b) We say that $u$ is $l$-$placed$ if $u$ is placed and $l(u)\! =\! l$. We say that $u$ is 
$(\leq\! l)$-$placed$ (resp., $(<\! l)$-$placed$, $(>\! l)$-$placed$) if there is $l'\!\leq\! l$ (resp., 
$l'\! <\! l$, $l'\! >\! l$) such that $u$ is $l'$-placed.\end{defi}

 When we consider the finite approximations of an element of $\mathbb{A}_3$, we have to guess the natural number $n$. We usually make some mistakes. In this case, we have to be able to come back to an earlier position. This is the role of the following predecessors.\bigskip

\noindent\bf Notation.\rm\ Let $u\!\in\! 2^{<\omega}$. Note that $<\eta >$ is $0$-placed with 
$\varepsilon (<\eta >)\! =\!\eta$ if $\eta\!\in\! 2$. This allows us to define
$$u^{-}\! :=\!\left\{\!\!\!\!\!\!\!
\begin{array}{ll}
& \emptyset\mbox{ if }\vert u\vert\!\leq\! 1\mbox{,}\cr
& u\vert\mbox{max}\{ l\! <\!\vert u\vert\mid u\vert l\mbox{ is placed}\}\mbox{ if }\vert u\vert\!\geq\! 2\mbox{,}
\end{array}
\right.$$
and, for $l\!\in\!\omega$,
$$u^{-l}\! :=\!\left\{\!\!\!\!\!\!\!
\begin{array}{ll}
& \emptyset\mbox{ if }\vert u\vert\!\leq\! 1\mbox{,}\cr
& u\vert\mbox{max}\{ k\! <\!\vert u\vert\mid u\vert k\mbox{ is }(\leq\! l)\mbox{-placed}\}\mbox{ if }
\vert u\vert\!\geq\! 2.
\end{array}
\right.$$
The following key lemma explains the relation between these predecessors and the placed sequences.

\begin{lem} Let $l\!\in\!\omega$ and $u\!\in\! 2^{<\omega}$ be $l$-placed with $\vert u\vert\!\geq\! 2$.\smallskip

(a) Assume that $u^-$ is $l$-placed. Then $\varepsilon (u^-)\! =\!\varepsilon (u)$. If moreover 
$(u^l)^-$ is $l$-placed, then the equality $(u^l)^-\! =\! (u^-)^l$ holds.\smallskip

(b) $u^{-l}$ is $l$-placed if and only if $(u^l)^{-l}$ is $l$-placed. In this case, 
$\varepsilon (u^{-l})\! =\!\varepsilon (u)$ and the equality $(u^l)^{-l}\! =\! (u^{-l})^l$ holds.\smallskip

(c) Assume that $u^-$ or $(u^l)^-$ is $(<\! l)$-placed. Then 
$u^-\! =\! u^{-l}\! =\! (u^l)^-\! =\! (u^l)^{-l}$.\smallskip

(d) Assume that $u^-$ or $(u^l)^-$ is $(>\! l)$-placed. Then exactly one of those two sequences is 
$(>\! l)$-placed, and the other one is $l$-placed. If $u^-$ (resp., $(u^l)^-$) is $(>\! l)$-placed, then 
$u^{-l}\! =\!\big( (u^l)^-\big)^l$ (resp., $u^{-l}\! =\! u^-$) and $\varepsilon (u^{-l})\! =\!\varepsilon (u)$ (resp., $\varepsilon\big( (u^l)^{-l}\big)\! =\!\varepsilon (u^l)$).\end{lem}

\noindent\bf Proof.\rm\ We first prove the following claim:\bigskip

\noindent\bf Claim.\it\ (i) Assume that $(\vert u^m\vert )_1\! =\! 0$. Then 
$u^-\! =\! u^{-l}\! =\! (u^l)^-\! =\! (u^l)^{-l}$ is $(<\! l)$-placed.\smallskip

(ii) Assume that $(\vert u^m\vert )_1\! >\! 0$. Then $u^-$ (resp., $u^{-l}$) is $(\geq\! l)$-placed and there is $j_0$ (resp., $j_1$) with $u^-\! =\! u\vert (<l(u^-),j_0>\! +1)$ (resp., $u^{-l}\! =\! u\vert (<l,j_1>\! +1)$).\bigskip

\noindent\bf Proof.\rm\ (i) Note that $l\!\geq\! 1$ since $\vert u\vert\!\geq\! 2$. As 
$(\vert u^m\vert )_1\! =\! 0$, $\vert u^m\vert\! =<(\vert u^m\vert )_0,(\vert u^m\vert )_1>=<l(u),0>$ and the sequence $u^-$ is $(<\! l)$-placed, which implies that $u^-\! =\! u^{-l}\! =\! (u^l)^-\! =\! (u^l)^{-l}$.\bigskip

\noindent (ii) The last assertion about $j_0$ and $j_1$ comes from the first one. It is enough to see that 
$u^-$ is $(\geq\! l)$-placed since the proof for $u^{-l}$ is similar. We argue by contradiction. Then 
$u\vert (<l,0>\! +1)$ is $l$-placed and 
$u\vert (<l,0>\! +1)\!\subsetneqq\! u\vert (<l,(\vert u^m\vert )_1>\! +1)\!\subseteq\! u$, so that 
$u\vert (<l,0>\! +1)\!\subsetneqq\! u^-$. This implies that 
$l\! +\! 0\!\leq\! l(u^-)\! +\! (\vert u^-\vert\! -\! 1)_1$, $(\vert u^-\vert\! -\! 1)_1\!\geq\! l\! -\! l(u^-)\! >\! 0$ 
and $u^-(\vert u^-\vert\! -\! 1)\! =\! 1$. But 
$$\begin{array}{ll}
u^-(\vert u^-\vert\! -\! 1)\!\!\!
& \! =\! u^-(<l(u^-),(\vert u^-\vert\! -\! 1)_1>)\! =\! u(<l(u^-),(\vert u^-\vert\! -\! 1)_1>)\cr
& \! =\! (u)_{l(u^-)}\big( (\vert u^-\vert\! -\! 1)_1\big)
\! =\! (s_{t_l( l(u^-))}10^\infty )\big( (\vert u^-\vert\! -\! 1)_1\big) .
\end{array}$$
Lemma 4.5.(a) implies that $(\vert u^-\vert\! -\! 1)_1\! <\! l\! -\! l(u^-)$, which is absurd.\hfill{$\diamond$}\bigskip

\noindent (a) By the claim, $(\vert u^m\vert )_1\! >\! 0$. Therefore 
$u\vert (<l,0>\! +1)\!\subsetneqq\! u\vert (<l,(\vert u^m\vert )_1>\! +1)\!\subseteq\! u$ is $l$-placed, 
$u\vert (<l,0>\! +1)\!\subseteq\! u^-$ and $<l,0> <\!\vert u^-\vert$. Thus 
$\varepsilon (u^-)\! =\! (u^-)(<l,0>)\! =\! u(<l,0>)\! =\!\varepsilon (u)$.\bigskip

 Assume now that $(u^l)^-$ is $l$-placed. As $u\vert (<l,0>\! +1)\!\subseteq\! u^-\!\subsetneqq\! u$, we get 
$$\big( u\vert (<l,0>\! +1)\big)^l\!\subseteq\! (u^l)^-.$$  
Thus $<l,0><\!\vert (u^l)^-\vert$. If $u^-\! =\! u\vert (<l,j_0>\! +1)$, then there is no 
$j_0\! <\! j\! <\! (\vert u^m\vert )_1$ such that $u(<l,j>)\! =\! 1$, and 
$(u^l)^-\! =\! u^l\vert (<l,j_0>\! +1)\! =\! (u^-)^l$.\bigskip

\noindent (b) Assume that $u^{-l}$ is $l$-placed. By the claim, we get $(\vert u^m\vert )_1\! >\! 0$ and $j_1$ with 
$$u^{-l}\! =\! u\vert (<l,j_1>\! +1).$$ 
Thus $(u^l)^{-l}\! =\! u^l\vert (<l,j_1>\! +1)\! =\! (u^{-l})^l$ is $l$-placed, by Lemma 4.8. The equivalence comes from the fact that $(u^l)^l\! =\! u$. We argue as in (a) to see that 
$\varepsilon (u^{-l})\! =\!\varepsilon (u)$ if $u^{-l}$ is $l$-placed.\bigskip

\noindent (c) Assume first that $u^-$ is $(<\! l)$-placed. Then $(\vert u^m\vert )_1\! =\! 0$, by the claim, (ii). Now the claim, (i), gives the result. If $(u^l)^-$  is $(<\! l)$-placed, then we apply this to $u^l$, using the facts that $u^l$ is $l$-placed and $(u^l)^l\! =\! u$.\bigskip

\noindent (d) Assume first that $u^-$ is $(>\! l)$-placed. The claim, (i), implies that 
$(\vert u^m\vert )_1\! >\! 0$, and the claim, (ii), gives $j_1$ with $u^{-l}\! =\! u\vert (<l,j_1>\! +1)$. Note that $u^{-l}\subsetneqq\! u^-$, $(u^-)_l\!\subseteq\! s_{t_{l(u^-)}(l)}10^\infty$ and 
$M(s_{t_{l(u^-)}(l)}10^\infty)\! <\! l(u^-)\! -\! l$, by Lemma 4.5.(a). Thus 
$$<l,M(s_{t_{l(u^-)}(l)}10^\infty )>\leq <l(u^-),0>\leq <l(u^-),(\vert u^-\vert\! -\! 1)_1>=
\!\vert u^-\vert\! -\! 1$$
and $(u^-)_l\big( M(s_{t_{l(u^-)}(l)}10^\infty )\big)$ is defined. This shows that 
$j_1\! =\! M( s_{t_{l(u^-)}(l)}10^\infty )$.\bigskip

 Note that $u^l\vert (<l,j_1>\! +1)\!\subseteq\! (u^l)^-$. The claim, (ii), shows that 
$(u^l)^{-l}\! =\! u^l\vert (<l,j_1>\! +1)$. We argue by contradiction to see that 
$(u^l)^-$ is not $(>\! l)$-placed. The proof of the previous point shows that 
$j_1\! =\! M( s_{t_{l( (u^l)^-)}(l)}10^\infty )$. Lemma 4.5.(b) shows that 
$s_{t_{l(u^-)}(l)}10^\infty\! =\! s_{t_{l( (u^l)^-)}(l)}10^\infty$. Thus 
$(u^-)_l(0)\! =\! (s_{t_{l(u^-)}(l)}10^\infty )(0)\! =\! (s_{t_{l((u^l)^-)}(l)}10^\infty )(0)\! =\!\big( (u^l)^-\big)_l(0)$, 
$$\varepsilon (u)\! =\! u(<l,0>)\! =\! (u)_l(0)\! =\! (u^-)_l(0)\! =\!\big( (u^l)^-\big)_l(0)\! =\!\varepsilon (u^l)\mbox{,}$$ 
which is absurd. This shows that $(u^l)^-\! =\! u^l\vert (<l,j_1>\! +1)\! =\! (u^{-l})^l$ is $l$-placed, by Lemma 4.8, so that $u^{-l}\! =\!\big( (u^l)^-\big)^l$. Moreover, 
$\varepsilon (u^{-l})\! =\! (u^{-l})(<l,0>)\! =\! u(<l,0>)\! =\!\varepsilon (u)$.\bigskip

 Assume now that $(u^l)^-$ is $(>\! l)$-placed. As $u^l$ is $l$-placed and $(u^l)^l\! =\! u$, the previous arguments show that $u^-$ is $l$-placed. In particular, $u^{-l}\! =\! u^-$.\hfill{$\square$}
 
\begin{thm} $(\mathbb{X}_3,\mathbb{A}_3)$ satisfies the conjecture.\end{thm}

\noindent\bf Proof.\rm\ We already noticed that it is enough to see that (a) or (b) holds. In Condition (5) in the proof of Theorem 5.1 in [L-Z], $u^{-l}$ should be replaced with $u^{-l(u)}$. We need to check that the map $f$ defined there satisfies $\mathbb{A}_3\!\subseteq\! (f\!\times\! f)^{-1}(A)$. So let 
$(\alpha ,\beta )\!\in\!\mathbb{A}_3$, which defines $n$. Let $(p_j)_{j\in\omega}$ be the infinite strictly increasing sequence of natural numbers $p_j\!\geq\! 1$ such that 
$(p_j\! -\! 1)_0\! =\! n$, $(p_j\! -\! 1)_1\! >\! 0$ and $\alpha (p_j\! -\! 1)\! =\! 1$. In particular, 
$\alpha\vert p_j$ is $n$-placed and $\varepsilon (\alpha\vert p_j)\! =\! 0$. Note that $(p_j)_{j\in\omega}$ is also the infinite strictly increasing sequence of natural numbers $p_j\!\geq\! 1$ such that 
$(p_j\! -\! 1)_0\! =\! n$, $(p_j\! -\! 1)_1\! >\! 0$ and $\beta (p_j\! -\! 1)\! =\! 1$ on one side, and a subsequence of both $(p^\alpha_k)_{k\in\omega}$ and $(p^\beta_k)_{k\in\omega}$ on the other side.\bigskip

 If moreover $p\!\geq\! p_0$ and $\alpha\vert p$ is placed, then $l(\alpha\vert p)\!\geq\! n$, by Lemma 4.8. In particular, if $p\!\geq\! p_0$ and $\alpha\vert p$ is $(\leq\! n)$-placed, then $\alpha\vert p$ is $n$-placed. This proves that $(p_j)_{j\in\omega}$ is the infinite strictly increasing sequence of integers 
$p_j\!\geq\! p_0$ such that $\alpha\vert p_j$ is $(\leq\! n)$-placed. Therefore 
$(\alpha\vert p_{j+1})^{-n}\! =\!\alpha\vert p_j$.\bigskip

 By Condition (3), $(U_{\alpha\vert p_j})_{j\in\omega}$ is a non-increasing sequence of nonempty clopen subsets of $A\cap\Omega_{X^2}$ whose GH-diameter tend to $0$. So we can define 
$F(\alpha ,\beta )\!\in\! A$ by $\{ F(\alpha ,\beta )\}\! :=\!\bigcap_{j\in\omega}~U_{\alpha\vert p_j}$. Note that 
$F(\alpha ,\beta )\! =\!\mbox{lim}_{j\rightarrow\infty}~(x_{\alpha\vert p_j},x_{\beta\vert p_j})\! =\!
\big( f(\alpha ),f(\beta )\big)\!\in\! A$, so that $\mathbb{A}_3\!\subseteq\! (f\!\times\! f)^{-1}(A)$.\bigskip
 
 It remains, when $k\!\geq\! 2$ (second case), to replace $l\! -\! 1$ with $l(u^-)$.\hfill{$\square$}\bigskip
 
\noindent\bf The general case\rm\bigskip
 
 Here we just give, for each $i\!\in\!\omega$, a sequence $(W_{\eta_i,n})_{n\in\omega}$ of nonempty 
$T_{\eta_i}$-open sets such that 
$Q_{\eta_i}\!\subseteq\!\overline{\bigcup_{n\in\omega}~W_{\eta_i,n}}^{T_{\eta_i}}$. This will imply that 
$\mathbb{K}_\xi$ has no $\borxi$-measurable countable coloring, by Corollary 4.3. We assume that 
$\xi\!\geq\! 4$, so that we may assume that $\eta_i\!\geq\! 3$. If 
$\eta\! =\!\mbox{sup}_{n\in\omega}~(\theta_n\! +\! 1)\!\geq\! 2$, then we set 
$V_{\eta ,n}\! :=\!\{\alpha\!\in\! 2^\omega\mid\forall i\! <\! n~~(\alpha )_i\!\notin\! P_{\theta_i}\ \wedge\ 
(\alpha )_n\!\in\! P_{\theta_n}\}$. We set, for $\eta\!\geq\! 3$, 
$$W_{\eta ,n}\! :=\!\big\{\alpha\!\in\! 2^\omega\mid\alpha (0)\! =\! s_{n+1}(0)\ \wedge\ 
(\alpha^*)_n\!\in\! P_{\theta_n}\ \wedge\ \forall i\! <\! n~~
(\alpha^*)_i\!\in\!\bigcup_{j<n-i}~V_{\theta_i,j}\big\} .$$
M\'atrai's construction ensures that $V_{\eta ,n}$ is $\tau_\eta$-open, and that $W_{\eta ,n}$ is a nonempty $T_\eta$-open set. Let $O$ be a basic $T_\eta$-open set meeting $Q_\eta$. As 
$T_\eta\! =\!\mbox{discrete}\!\times\!\tau_\eta$ and 
${\tau_\eta}_{\vert P_\eta}\!\equiv\! (\Pi_{i\in\omega}~{\tau_{\theta_i}})_{\vert P_\eta}$, we can find 
$\varepsilon\!\in\! 2$ and $(O_i)_{i<l}\!\in\!\pi_{i<l}~\tau_{\theta_i}$ such that 
$O\! =\!\{\alpha\!\in\! 2^\omega\mid
\alpha (0)\! =\!\varepsilon\ \wedge\ \forall i\! <\! l~~(\alpha^*)_i\!\in\! O_i\}$. As $P_{\theta_i}$ is 
$\tau_{\theta_i}$-closed nowhere dense and 
$\neg P_{\theta_i}\! =\!\bigcup_{n\in\omega}~V_{\theta_i,n}$, we can find $n_i$ such that $O_i$ meets $V_{\theta_i,n_i}$. We choose $n\! >\!\mbox{max}_{i<l}~(n_i\! +\! i)$ such that 
$s_{n+1}(0)\! =\!\varepsilon$. Then $W_{\eta ,n}$ meets $O$.\bigskip

 Our motivation to introduce these examples is that they induce a set $\mathbb{K}_3$ satisfying the conjecture. This is the reason why we think that they are reasonable candidates for the general case.\bigskip
 
\noindent\bf (B) The small classes\rm\bigskip
 
 In Section 3, we met $D_2(\bormone )$ graphs of fixed point free partial injections with a Borel countable (2-)coloring, but without $\borxi$-measurable countable coloring. Their complement are 
$\check D_2(\bormone )$ sets in $(\borel\!\times\!\boraone )_\sigma$, but not in 
$(\boraxi\!\times\!\boraxi )_\sigma$. However, a positive result holds for the simpler classes, which shows some optimality in our results.

\begin{prop} Let ${\bf\Gamma}\!\subseteq\! D_2(\bormone)$ be a Wadge class (in zero-dimensional spaces), and $A$ be a set in ${\bf\Gamma}\cap (\borel\!\times\!\boraone )_\sigma$ (resp., $(\borel\!\times\!\borel )_\sigma$). Then $A\!\in\! ({\bf\Gamma}\!\times\!\boraone )_\sigma$ (resp., $({\bf\Gamma}\!\times\! {\bf\Gamma})_\sigma$).\end{prop}
 
\noindent\bf Proof.\rm\ Let us do it for $(\borel\!\times\!\boraone )_\sigma$, the other case being similar. The result is clear for $\{\emptyset\}$, $\check\{\emptyset\}$, $\borone$, 
$\boraone$. If ${\bf\Gamma}\! =\!\bormone$, then we can write 
$A\! =\!\bigcup_{n\in\omega}~C_n\!\times\! D_n$, with $C_n\!\in\!\borel$ and $D_n\!\in\!\boraone$. We just have to note that $A\! =\!\bigcup_{n\in\omega}~\overline{C_n}\!\times\! D_n$. If 
${\bf\Gamma}\! =\!\bormone\oplus\boraone$, then we can write 
$A\! =\!\bigcup_{n\in\omega}~C_n\!\times\! D_n\! =\! (C\cap D)\cup (O\!\setminus\! D)$, with 
$C_n\!\in\!\borel$, $\neg C, O, D,\neg D, D_n\!\in\!\boraone$. Note that 
$A\! =\! (D\cap\bigcup_{n\in\omega}~\overline{C_n}\!\times\! D_n)\cup (O\!\setminus\! D)$. Finally, if 
${\bf\Gamma}\! =\! D_2(\bormone )$, then write 
$A\! =\!\bigcup_{n\in\omega}~C_n\!\times\! D_n\! =\! C\cap O$, with 
$C_n\!\in\!\borel$, $\neg C, O, D_n\!\in\!\boraone$. Note that 
$A\! =\! O\cap\bigcup_{n\in\omega}~\overline{C_n}\!\times\! D_n$.\hfill{$\square$}\bigskip

\noindent\bf (C) The finite case\rm
 
\begin{prop} Assume that $\bf\Gamma$ is closed under finite intersections and continuous pre-images, $X,Y$ are topological spaces, $\kappa$ is finite, and $A\!\in\! {\bf\Gamma}(X\!\times\! Y)$ is the union of $\kappa$ rectangles. Then $A$ is the union of at most $2^{2^\kappa}$ rectangles whose sides are in $\bf\Gamma$.\end{prop}

\noindent\bf Proof.\rm\ Assume that $A\! =\!\bigcup_{n<\kappa}~A_n\!\times\! B_n$. Let us prove that 
$$A\! =\!\bigcup_{I\subseteq\kappa, (\bigcap_{n\in I}~A_n)\setminus (\bigcap_{n\notin I}~A_n)
\not=\emptyset}~(\bigcap_{n\in I}~A_n)\!\times\! (\bigcup_{n\in I}~B_n).$$ 
So let $(x,y)\!\in\! A$, and let $I\! :=\!\{ n\! <\!\kappa\mid x\!\in\! A_n\}$. Then 
$x\!\in\! (\bigcap_{n\in I}~A_n)\!\setminus\! (\bigcap_{n\notin I}~A_n)$, and  
$(x,y)$ is in $(\bigcap_{n\in I}~A_n)\!\times\! (\bigcup_{n\in I}~B_n)$ since 
$(x,y)\!\in\! A_n\!\times\! B_n$ for some $n\! <\!\kappa$. The other inclusion is clear.

\vfill\eject

 Assume now that $x\!\in\! (\bigcap_{n\in I}~A_n)\!\setminus\! (\bigcap_{n\notin I}~A_n)$. Then 
$\bigcup_{n\in I}~B_n\! =\! A_x\! =\! f^{-1}(A)$, where the formula $f(y)\! :=\! (x,y)$ defines 
$f\! :\! Y\!\rightarrow\! X\!\times\! Y$ continuous. This shows that $\bigcup_{n\in I}~B_n$ is in $\bf\Gamma$. So we proved the following:\bigskip

 $A$ is the union of at most $2^\kappa$ rectangles $A'_n\!\times\! B'_n$, where $A'_n$ is a finite intersection of some of the $A_n$'s, and $B'_n$ is a finite union of some of the $B_n$'s which is in $\bf\Gamma$.\bigskip
 
 Applying this again, we see that $A$ is the union of at most $2^{2^\kappa}$ rectangles 
$A''_n\!\times\! B''_n$, where $A''_n$ is a finite union of some of the $A'_n$'s which is in $\bf\Gamma$, and $B''_n$ is a finite intersection of some of the $B'_n$'s. We are done since $\bf\Gamma$ is closed under finite intersections.\hfill{$\square$}\bigskip

 This proof also shows the following result:

\begin{prop} Assume that $\bf\Gamma$ is closed under continuous pre-images, $X,Y$ are topological spa-ces, $\kappa$ is finite, and $A\!\in\! {\bf\Gamma}(X\!\times\! Y)$ is the union of 
$\kappa$ rectangles of the form $2^X\!\times\!\boraone (Y)$. Then $A$ is the union of at most 
$2^{2^\kappa}$ rectangles of the form ${\bf\Gamma}(X)\!\times\!\boraone (Y)$.\end{prop}

\noindent\bf Remarks.\rm\ (1) For colorings, Theorem 1.2 gives, for each $\xi$, a 
$D_2(\bormone )$ binary relation with a Borel finite (2-)coloring, but with no $\borxi$-measurable finite coloring.\bigskip

\noindent (2) $\emptyset$ has a 1-coloring. An open binary relation having a finite coloring $c$ has also a $D_2(\bormone )$-measurable finite coloring (consider the differences of the 
$\overline{c^{-1}(\{ n\} )}$'s, for $n$ in the range of the coloring). This leads to the following question:\bigskip 
 
\noindent\bf Question.\rm\ Can we build, for each $\xi$, a closed binary relation with a Borel finite coloring but no $\borxi$-measurable finite coloring?

\section{$\!\!\!\!\!\!$ References}

\noindent [K]\ \ A. S. Kechris,~\it Classical Descriptive Set Theory,~\rm 
Springer-Verlag, 1995

\noindent [K-S-T]\ \ A. S. Kechris, S. Solecki and S. Todor\v cevi\'c, Borel chromatic numbers,\ \it 
Adv. Math.\rm\ 141 (1999), 1-44

\noindent [L-Z]\ \ D. Lecomte and M. Zelen\'y, Baire-class $\xi$ colorings: the first three levels,\ \it to appear in Trans. Amer. Math. Soc. (see arXiv)\rm\ 

\noindent [Lo]\ \ A. Louveau, A separation theorem for $\Ana$ sets,\ \it Trans. 
Amer. Math. Soc.\ \rm 260 (1980), 363-378

\noindent [M\'a]\ \ T. M\'atrai, On the closure of Baire classes under transfinite convergences,~\it Fund. Math.\ \rm 183, 2 (2004), 157-168
 
\end{document}